\documentclass{memoir} 

\setlrmarginsandblock{1.75in}{1.75in}{*}
\usepackage{amsmath,amstext,amssymb,amscd}

\numberwithin{figure}{section}

\usepackage{hyperref}
\usepackage[hyperref, amsmath, amsthm,thmmarks]{ntheorem}
\usepackage[hyperref, ntheorem]{cleveref}

\usepackage{graphicx}
\usepackage{relsize}

\tightlists 

\renewcommand{\.}[1]{\hbox to 0pt{#1\hss}}

\AtBeginDocument{%
\setlength{\abovedisplayskip}{0.25em plus 0.1em}
\setlength{\belowdisplayskip}{0.25em plus 0.1em}
}

\numberwithin{equation}{section}

\makeatletter
\renewtheoremstyle{plain}%
{\item[\hskip\labelsep \normalfont(##2) \theorem@headerfont ##1\@addpunct{.}]}%
{\item[\hskip\labelsep \normalfont(##2) \theorem@headerfont ##1 \normalfont ##3\@addpunct{.}]}

\renewtheoremstyle{definition}%
{\item[\hskip\labelsep \normalfont(##2) \theorem@headerfont ##1\@addpunct{.}]}%
{\item[\hskip\labelsep \normalfont(##2) \theorem@headerfont ##1 \normalfont ##3\@addpunct{.}]}
\makeatother

\crefformat{section}{Section~#2#1#3}
\crefformat{subsection}{\S#2#1#3}
\crefformat{figure}{Figure~#2#1#3}
\crefmultiformat{section}{Sections~#2#1#3}{, #2#1#3}{ and~#2#1#3}

\theoremstyle{plain}

\newtheorem{thm}[equation]{Theorem}

\newtheorem{lem}[equation]{Lemma}
\newtheorem{cor}[equation]{Corollary}
\newtheorem{prop}[equation]{Proposition}

\newtheorem{conjecture}[equation]{Conjecture}

\crefformat{equation}{Equation~(#2#1#3)}
\crefformat{thm}{Theorem~#2#1#3}
\crefformat{mainthm}{Theorem~#2#1#3}
\crefformat{lem}{Lemma~#2#1#3}
\crefformat{cor}{Corollary~#2#1#3}
\crefformat{prop}{Proposition~#2#1#3}
\crefformat{openproblem}{Problem~#2#1#3}
\crefformat{egprop}{Example~#2#1#3}
\crefformat{conjecture}{Conjecture~#2#1#3}
\crefformat{keyfact}{Key Fact~#2#1#3}

\crefformat{exer}{Exercise~#2#1#3} 
\crefformat{exeri}{Exercise~#2#1#3} 
\crefformat{exerii}{Exercise~#2#1#3}
\crefformat{exeriii}{Exercise~#2#1#3}

\crefrangeformat{exeri}{Exercises~(#3#1#4)--(#5#2#6)}

\crefmultiformat{exeri}{Exercises~#2#1#3}{, #2#1#3}{ and~#2#1#3}
\crefmultiformat{rem}{Remarks~#2#1#3}{, #2#1#3}{ and~#2#1#3}
\crefmultiformat{cor}{Corollaries~#2#1#3}{, #2#1#3}{ and~#2#1#3}
\crefmultiformat{thm}{Theorems~#2#1#3}{, #2#1#3}{ and~#2#1#3}
\crefmultiformat{lem}{Lemmas~#2#1#3}{, #2#1#3}{ and~#2#1#3}
\crefmultiformat{prop}{Propositions~#2#1#3}{, #2#1#3}{ and~#2#1#3}
\crefmultiformat{conjecture}{Conjectures~#2#1#3}{, #2#1#3}{ and~#2#1#3}

\crefformat{MostRigThm}{Theorem~#2#1#3}
\crefformat{MargSuperThm}{Theorem~#2#1#3}
\crefformat{MargArithThm}{Theorem~#2#1#3}

\theoremstyle{definition}
\theorembodyfont{\normalfont}

\newtheorem{defn}[equation]{Definition}
\newtheorem{notation}[equation]{Notation}
\newtheorem{warn}[equation]{Warning}

\newtheorem{eg}[equation]{Example}

\newtheorem{assump}[equation]{Assumption}
\newtheorem{problem}[equation]{Problem}

\crefformat{defn}{Definition~#2#1#3}
\crefformat{notation}{Notation~#2#1#3}
\crefformat{terminology}{Terminology~#2#1#3}
\crefformat{exmp}{Example~#2#1#3}
\crefformat{eg}{Example~#2#1#3}
\crefformat{prob}{Problem~#2#1#3}
\crefformat{assump}{Assumption~#2#1#3}
\crefformat{problem}{Problem~#2#1#3}
\crefformat{warn}{Warning~#2#1#3}

\crefmultiformat{eg}{Examples~#2#1#3}{, #2#1#3}{ and~#2#1#3}
\crefmultiformat{defn}{Definitions~#2#1#3}{, #2#1#3}{ and~#2#1#3}

\theoremstyle{remark}

\newtheorem{rem}[equation]{Remark}       
\newtheorem{remgen}[equation]{Remark on the general case}

\crefformat{rem}{Remark~#2#1#3}
\crefformat{remgen}{Remark~#2#1#3}

\newtheorem{exer}[equation]{Exercise}

 \newcounter{case}
 \newenvironment{case}[1][\unskip]{\refstepcounter{case}\em
 \medskip \noindent Case \thecase\ #1.\ }{\unskip\upshape}
 \renewcommand{\thecase}{\arabic{case}}

\crefformat{case}{Case~#2#1#3}

 \newcounter{subcase}
 
 \numberwithin{subcase}{case}

\crefformat{subcase}{Subcase~#2#1#3}

 \newcounter{step}
 
 \renewcommand{\thestep}{\arabic{step}}

\crefformat{step}{Step~#2#1#3}
\crefmultiformat{step}{Steps~#2#1#3}{, #2#1#3}{ and~#2#1#3}

\makeatletter
\newenvironment{thmref}{\thmrefer}{}
\newcommand{\thmrefer}[1]{\expandafter\ifx\csname r@#1\endcsname\relax
 \renewcommand{\theequation}{??}
     \protect\G@refundefinedtrue%
    \nfss@text{\reset@font\bfseries ??}%
    \@latex@warning{Reference `#1' on page \thepage \space undefined}%
 \else \cref@getlabel{#1}{\dtemp} \renewcommand{\theequation}{\dtemp$'$} \fi
  \addtocounter{equation}{-1}
}
 \makeatother

\newcommand{\thmrefrefer}[1]{\renewcommand\theequation
  {\protect\ref{#1}$''$}\addtocounter{equation}{-1}}

\renewcommand{\pref}[1]{{\upshape(\ref{#1})}} 

\newcommand{\fullref}[2]{\ref{#1}\pref{#1-#2}}
\newcommand{\fullcref}[2]{\cref{#1}\pref{#1-#2}}

\DeclareMathOperator{\Stab}{Stab}
\DeclareMathOperator{\supp}{supp}

\newcommand{\Isom}{\operatorname{Isom}\nolimits}

\newcommand{\fund}{\mathord{\mathcal{F}}}
\newcommand{\flags}{\mathord{\mathcal{F}}}

\newcommand{\TVS}{\mathop{\mathcal{V}}\nolimits}

\newcommand{\Cbdd}{C_{\text{\normalfont\upshape{bdd}}}}
\newcommand{\hCbdd}{\dot C_{\text{\normalfont\upshape{bdd}}}}
\newcommand{\Hbdd}{H_{\text{\normalfont\upshape{bdd}}}}

\renewcommand{\operatorname}[1]{\mathop{\mathrm{#1}}\nolimits}

\newcommand{\Prob}{\operatorname{Prob}}

\newcommand{\GL}{\operatorname{GL}}

\newcommand{\Id}{\operatorname{Id}}

\newcommand{\SL}{\operatorname{SL}}

\newcommand{\SO}{\operatorname{SO}}
\newcommand{\Sp}{\operatorname{Sp}}

 \newcommand{\Func}{\mathop{\hbox{$\mathcal{F}$}}\nolimits}

\newcommand{\real}{{\mathord{\mathbb R}}}
\newcommand{\rational}{{\mathord{\mathbb Q}}}
\newcommand{\complex}{{\mathord{\mathbb C}}}

\newcommand{\integer}{{\mathord{\mathbb Z}}}

\renewcommand{\natural}{\mathord{\mathbb N}}
\renewcommand{\circle}{S^1}

\DeclareMathOperator{\Homeo}{Homeo}
\DeclareMathOperator{\Diff}{Diff}

\newcommand{\Qrank}
 {\mathop{\mathbb{Q}\text{\upshape{-rank}}}\nolimits}

\newcommand{\Rrank}
 {\mathop{\mathbb{R}\text{\upshape{-rank}}}\nolimits}

\newcommand\bigset[2]{\left\{\, #1 
 \mathrel{\left| \vphantom {\left\{ #1 \mid #2 \right\} }
 \right.} #2 \,\right\} }

\newcommand{\hilbert}{\mathord{\mathcal{H}}}

\newcommand{\iso}{\cong}

\newcommand{\precc}{\mathrel{\prec\mkern-8mu\prec}}

\newcommand{\defit}{\emph}

\renewcommand{\see}[1]{(see~\ref{#1})}
\newcommand{\csee}[1]{(see~\cref{#1})}

\newcommand{\cf}[1]{(cf.~\ref{#1})}
\newcommand{\ccf}[1]{(cf.\ \cref{#1})}

\newcommand{\fullcsee}[2]{{\upshape(}see \cref{#1}(\ref{#1-#2}){\upshape)}}

\newcommand{\hint}[1]{\par\noindent{\smaller{[}\emph{Hint:} #1{]}\par}}

\renewcommand{\thesection}{\arabic{section}}

\renewcommand{\thesubsection}{\thesection\Alph{subsection}}

\makeatletter
\newenvironment{notes}
 {\medbreak\@startsection{notesection}{3}
  \z@\z@{-.5em}%
  {\bfseries}*{Notes for \S\thesection. }}{}

\newenvironment{subsecnotes}
 {\medbreak\@startsection{notesection}{3}
  \z@\z@{-.5em}%
  {\bfseries}*{Notes for \S\thesubsection. }}{}
\makeatother

\begin{document}

\title{Can Lattices in $\SL(n,\real)$ Act on the Circle?}

\author{Dave Witte Morris}

\renewcommand{\maketitlehookd}{\begin{center} \itshape To Robert J.\ Zimmer on his 60th birthday\end{center}}

\maketitle

\section{Overview} \label{IntroSect}

The following theorem is a simple example of the Zimmer program's principle that large groups should not be able to act on small manifolds. (For a description of the Zimmer program, see D.\,Fisher's survey article \cite{Fisher-Survey}
elsewhere in this volume.) Unless stated otherwise, we assume actions are continuous, but we require no additional regularity.

\begin{thm}[(Witte)] \label{Witte-AOCThm}
Let\/ $\Gamma$ be a finite-index subgroup of\/ $\SL(n,\integer)$, with $n \ge 3$. Then\/ $\Gamma$ has no faithful action on the circle~$\circle$.
\end{thm}

\begin{rem} \ 
\noprelistbreak
\begin{enumerate}
\item The proof of this theorem is not at all difficult, and will be given in \cref{SLnZNotLOSect}, after the result is translated to a more algebraic form in \cref{AlgReformSect}. The proof illustrates the use of calculations with unipotent elements, which is a standard technique in the theory of arithmetic groups.
\item The assumption that $n \ge 3$ is essential. Indeed, some finite-index subgroup of $\SL(2,\integer)$ is a free group, which has countless actions on the circle (some faithful and some not).
\end{enumerate}
\end{rem}

A group~$\Gamma$ as in \cref{Witte-AOCThm} is a lattice in $\SL(n,\real)$. It is an open question whether the theorem generalizes to the other lattices:

\begin{conjecture} \label{LattSLnRNoActConj}
Let\/ $\Gamma$ be a lattice in\/ $\SL(n,\real)$, with $n \ge 3$. Then\/ $\Gamma$ has no faithful action on~$\circle$.
\end{conjecture}

Note that if this conjecture is true, then every action of\/~$\Gamma$ on $\circle$ or\/~$\real$ has a nontrivial kernel. Since the kernel is a normal subgroup, and the Margulis Normal Subgroup Theorem tells us that lattices in $\SL(n,\real)$ are ``almost simple'' (if $n \ge 3$), it follows that the kernel is a finite-index subgroup of~$\Gamma$. Therefore, every $\Gamma$-orbit is finite.
Such reasoning shows that the conjecture can be restated as follows: 

\begin{thmref}{LattSLnRNoActConj}
\begin{conjecture} \label{LattSLnRFinOrbConj}
Let\/ $\Gamma$ be a lattice in\/ $\SL(n,\real)$, with $n \ge 3$.  Whenever\/~$\Gamma$ acts on~$\circle$, every orbit is finite.
\end{conjecture}
\end{thmref}

There is considerable evidence for the above conjectures, including the following important theorem:

\begin{thm}[(Ghys, Burger-Monod)] \label{GhysBurgerMonodFinOrb}
If\/ $\Gamma$ is any lattice in\/ $\SL(n,\real)$, with $n \ge 3$, then every action of\/~$\Gamma$ on~$\circle$ has at least one finite orbit.
\end{thm}

\begin{rem}
The proof of \'E.\,Ghys utilizes amenability (or the Furstenberg boundary); see \cref{GhysPfSect,GhysMissingSect}. The proof of M.\,Burger and N.\,Monod is based on a vanishing theorem for bounded cohomology; see \cref{BurgerMonodPf}. 
\end{rem} 

The above conjectures are about continuous actions, but they can be weakened by considering only differentiable actions.
It will be explained in \cref{ReebThurstonSect} that combining the Ghys-Burger-Monod Theorem with the well-known Reeb-Thurston Stability Theorem establishes these weaker conjectures:

\begin{cor}[(Ghys, Burger-Monod)] \label{LattSLnNoC1Act}
If\/ $\Gamma$ is any lattice in\/ $\SL(n,\real)$, with $n \ge 3$, then:
	\noprelistbreak
	\begin{enumerate}
	\item $\Gamma$ has no faithful $C^1$ action on~$\circle$.
	\item Whenever\/ $\Gamma$ acts on~$\circle$ by $C^1$ diffeomorphisms, every orbit is finite.
	\end{enumerate}
\end{cor}

For $n \ge 3$, it is well known that every lattice in $\SL(n,\real)$ has Kazhdan's property $(T)$ \csee{KazhdanDefn}, and the following result shows that if we strengthen the
differentiability hypothesis slightly, then the conclusion is true
for all groups that have that property.  This is a very broad
class of groups (including many groups that are not even linear), so
one would expect a much stronger result to hold for the special case
of lattices in $\SL(n,\real)$. Thus, this theorem constitutes good evidence for
\cref{LattSLnRNoActConj}. The proof will be presented in \cref{NavasPfSect}; it is both elegant and elementary.

\begin{thm}[(Navas)] \label{Navas-AOC}
If\/ $\Gamma$ is any infinite, discrete group with Kazhdan's property~$(T)$, then\/ $\Gamma$ has no faithful $C^2$ action on the circle.
 \end{thm}

Bounded generation provides another approach to proving \cref{Witte-AOCThm}, and some other cases of \cref{LattSLnRNoActConj}. This strategy will be explained in \cref{BddOrbBddGenSect}.

In spite of the above results, \cref{LattSLnRNoActConj} remains completely open for cocompact lattices:

\begin{problem}
 Find a cocompact lattice~$\Gamma$ in $\SL(n,\real)$, for some~$n$, such that no finite-index
subgroup of~$\Gamma$ has a faithful action on~$\circle$.
 \end{problem}

\begin{rem}
A final section of the paper (\S\ref{ComplementsSect}) briefly discusses the generalization of \cref{LattSLnRNoActConj} to lattices in other semisimple Lie groups, and two other topics: 
lattices that do act on the circle, and actions on trees.
\end{rem}

\begin{notes}

The survey of \'E.~Ghys \cite{GhysCircleSurvey} and the forthcoming book of A.\,Navas \cite{Navas-GrpsDiffeos} are excellent
introductions to the study of group actions on the circle.
Versions of \cref{LattSLnRNoActConj} were discussed in conversation
as early as 1990, but the first published appearance may be in
\cite[p.~200]{GhysCircle}. 

Expositions of the Margulis Normal Subgroup Theorem can be found in \cite[Chap.~4]{MargulisBook} and \cite[Chap.~8]{ZimmerBook}.

Regarding the equivalence of \cref{LattSLnRNoActConj,LattSLnRFinOrbConj}, see \cite[Thm.~1]{UominiCpctLeafConj} for a proof that if every $\Gamma$-orbit on a connected manifold is finite, then the kernel of the action has finite index.

\end{notes}

\subsection*{Acknowledgments}
Preparation of this paper was partially supported by research grants from the Natural Sciences and Engineering Research Council of Canada and the U.~S.\ National Science Foundation.
I am grateful to 
M.\,Burger, \'E.\,Ghys, N.\,Monod, A.\,Navas, and an anonymous referee,
for helpful comments on a preliminary version of the manuscript.

\section{Algebraic reformulation of the conjecture} \label{AlgReformSect}

In this section, we explain that \cref{LattSLnRNoActConj} can be
reformulated in a purely algebraic form \see{AOC<>LO}.
 The proof has two main parts: having an action
on~$\circle$ is almost the same as having an action on~$\real$, and
having an action on~$\real$ is essentially the same as being left
orderable.

\begin{defn} \label{LODefn}
 $\Gamma$ is \defit{left orderable} if there is a left-invariant
total order on~$\Gamma$. More precisely, there is an order
relation~$\prec$ on~$\Gamma$, such that:
 \noprelistbreak
 \begin{enumerate}
 \item \label{LODefn-total}
 $\prec$ is a total order (that is, for all $a,b \in \Gamma$,
either $a \prec b$ or $a \succ b$, or $a =
b$);
 and
 \item \label{LODefn-mult}
 $\prec$ is invariant under multiplication on the left (that is, for all $a, b, c \in
\Gamma$, if $a \prec b$, then $c a \prec c b$).
 \end{enumerate}
 \end{defn}

\begin{exer} \label{LO<>PosCone}
 Show $\Gamma$ is left orderable if and only if there exists a subset~$P$
of~$\Gamma$, such that
 \noprelistbreak
 \begin{enumerate}
 \item $ \Gamma = P \sqcup \{e\} \sqcup P^{-1}$ (disjoint union),
 where $P^{-1} = \{\, a^{-1} \mid a \in P \,\}$; 
 and
 \item $P$ is closed under multiplication (that is, $a b
\in P$, for all $a, b \in P$).
 \end{enumerate}
 Thus, being left orderable is a property of the internal
algebraic structure of~$\Gamma$.
 \hint{A subset~$P$ as above is  the ``positive cone'' of an order~$\prec$. Given~$P$, define $a \prec b$ if $b^{-1} a
\in P$. 
Given~$\prec$, define $P = \{\, a \in \Gamma \mid a \succ e \,\}$.}
 \end{exer}

\begin{rem}
 The group~$\Gamma$ is said to be \defit{right orderable} if there is a
\emph{right}-invariant total order on~$\Gamma$. 
 The following exercise shows that the choice between ``left orderable" and ``right orderable" is
entirely a matter of personal preference.
 \end{rem}

\begin{exer} 
 Show that $\Gamma$ is left orderable if and only if $\Gamma$~is
right orderable.
 \hint{Define $x \precc y$ if $x^{-1} \prec y^{-1}$. Alternatively, the
conclusion can be derived from \cref{LO<>PosCone} and its
analogue for right orderable groups.}
\end{exer}

With this terminology, we can state the following algebraic conjecture, which will be seen to be equivalent to \cref{LattSLnRNoActConj}.

 \begin{conjecture} \label{GammaNotLO}
 Let\/ $\Gamma$ be a lattice in\/ $\SL(n,\real)$, with $n \ge 3$. Then\/
 $\Gamma$ is not left orderable.
 \end{conjecture}

As a tool for showing that \cref{GammaNotLO} is equivalent to \cref{LattSLnRNoActConj}, let us give a geometric interpretation of being left orderable.

\begin{lem} \label{LO<>AOR}
A countable group\/~$\Gamma$ is left orderable if and only if there is an orientation-preserving, faithful action of\/~$\Gamma$ on\/~$\real$.
 \end{lem}

\newcommand{\LOimpliesAOR}{\Rightarrow}
\newcommand{\AORimpliesLO}{\Leftarrow}

\begin{proof}
 ($\Leftarrow$) For $a,b \in \Gamma$, define
 $$ \mbox{$a \prec b$ \quad  if \quad $a(0) < b(0)$} .$$
 It is easy to see that $\prec$ is transitive and antisymmetric, so it
defines a partial order on~$\Gamma$.

 For any $c \in \Gamma$, the function $x \mapsto c(x)$ is a strictly
increasing function on~$\real$ (because the action is orientation
preserving). Hence, if $a(0) < b(0)$, then $c \bigl( a(0) \bigr) < c \bigl( b(0) \bigr)$. Therefore, $\prec$ is left-invariant (as
required in \fullref{LODefn}{mult}). 

However, the relation $\prec$ may not be a total order, because it could
happen that $a(0) = b(0)$, even though $a \neq
b$. It is not difficult to modify the definition to deal
with this technical point \csee{LO<>AOR-BreakTie}.

($\Rightarrow$) We leave this as an exercise for the reader\..%
 \footnote{\emph{Hint:} If the order relation $(\Gamma,\prec)$ is dense (that is, if
$a \prec b \Rightarrow \exists c \in \Gamma, a \prec c \prec
b$), then it is order-isomorphic to $(\rational,<)$. The action
of~$\Gamma$ on $(\Gamma,\prec)$ by left multiplication extends to a
continuous action of~$\Gamma$ on the Dedekind completion, which is
homeomorphic to~$\real$.}
 \end{proof}
 
 \begin{exer} \label{LO<>AOR-BreakTie}
 Complete the proof of \cref{LO<>AOR}($\AORimpliesLO$), by
defining a left-invariant total order on~$\Gamma$.
 \hint{Make a list $q_1,q_2,\ldots$ of the rational numbers, and
define $a \prec b$ if $a(q_n) < b(q_n)$, where
$n$~is minimal with $a(q_n) \neq b(q_n)$.}
\end{exer}

\begin{prop} \label{AOC<>LO}
 \Cref{LattSLnRNoActConj,GammaNotLO} are equivalent.
 \end{prop}

\begin{proof}
We show there exists a counterexample to \cref{LattSLnRNoActConj} if and only if there exists a counterexample to \cref{GammaNotLO}.

 ($\Leftarrow$) Suppose $\Gamma$ is a counterexample to \cref{GammaNotLO}, so $\Gamma$ is left orderable. Then \cref{LO<>AOR} tells us that $\Gamma$~has
a faithful action on~$\real$. This implies that
$\Gamma$ acts faithfully on the one-point compactification
of~$\real$, which is homeomorphic to~$\circle$. 
So $\Gamma$ is a counterexample to \cref{LattSLnRNoActConj}.

($\Rightarrow$) Suppose $\Gamma$ is a counterexample to \cref{LattSLnRNoActConj}, so  $\Gamma$ has a faithful action on~$\circle$.
From the Ghys-Burger-Monod Theorem \pref{GhysBurgerMonodFinOrb}, we know that the action has a finite orbit. Therefore, some subgroup~$\Gamma'$ of finite index in~$\Gamma$ has a fixed point~$p$. Then $\Gamma'$
acts faithfully on the complement $\circle \smallsetminus \{p\}$,
which is homeomorphic to~$\real$. Let $\Gamma''$ be the subgroup
of~$\Gamma'$ consisting of the elements that act by
orientation-preserving homeomorphisms (so $\Gamma''$ is either~$\Gamma'$, or a subgroup of index two in~$\Gamma'$). Then \cref{LO<>AOR} tells us that $\Gamma''$ is left orderable.
Furthermore, $\Gamma''$~is a lattice in $\SL(n,\real)$ (because it has finite index in~$\Gamma$).
Therefore, $\Gamma''$ is a counterexample to \cref{GammaNotLO}.
 \end{proof}

\begin{warn}
 The proof does not show that the assertions of the two conjectures are equivalent for
each individual lattice~$\Gamma$. Rather, if one of the conjectures is
valid for all lattices~$\Gamma$, then
the other conjecture is also valid for all~$\Gamma$.
 More precisely, if one of the conjectures is valid for every
finite-index subgroup of~$\Gamma$, then the other conjecture is also
valid for~$\Gamma$.
 \end{warn}

\begin{rem}
 In order to prove that \cref{GammaNotLO} implies
\cref{LattSLnRNoActConj}, we appealed to
\cref{GhysBurgerMonodFinOrb}, which is quite deep. A stronger algebraic
conjecture (that no central extension of~$\Gamma$ is left orderable)
can easily be shown to imply \cref{LattSLnRNoActConj}:
 \end{rem}
 
 \begin{exer} \label{>CentExtLO}
 A \defit{central extension} of~$\Gamma$ is a group~$\Lambda$, such
that
  $\Lambda / Z \iso \Gamma$, 
 for some subgroup~$Z$ of the center of~$\Lambda$.
Show that if $\Gamma$ has an orientation-preserving, faithful action on~$\circle$, then some central extension of~$\Gamma$ is left
orderable.
 \hint{Because $\real$ is the universal cover of~$\circle$, every
homeomorphism of~$\circle$ lifts to a homeomorphism of~$\real$, and
the lift is unique modulo an element of the fundamental
group~$\integer$. Let $\Lambda$ be the subgroup of $\Homeo_+(\real)$
consisting of all of the possible lifts of all of the elements
of~$\Gamma$, and show that $\Lambda$ is a central extension that is
left orderable.}
 \end{exer}

%





\begin{notes}
The material in this section is well known.
See \cite{KopytovMedvedev} for a treatment of the algebraic theory of left-ordered groups.
Informative discussions of orderings, actions on the line, and related topics appear in \cite[\S6.5]{GhysCircleSurvey} and \cite[\S2.2.3--2.2.6]{Navas-GrpsDiffeos}.
\end{notes}

\section{\texorpdfstring{$\SL(\text{\lowercase{$n$}},\integer)$}{SL(n,Z)} has no faithful action on the circle} \label{SLnZNotLOSect}

In this section, we prove \cref{Witte-AOCThm}, by exploiting the interaction between some obvious nilpotent subgroups of~$\Gamma$. The crucial ingredients are \Cref{LOHeis}, and the fact that the theorem can be restated in the following
algebraic form \cf{AOC<>LO}.

\begin{thm}[(Witte)] \label{SLnZnotLO}
 If\/ $\Gamma$ is a finite-index subgroup of\/ $\SL(n,\integer)$, with $n \ge 3$, then\/
$\Gamma$ is not left orderable.
 \end{thm}
 
 \begin{rem}
 It is easy to see that the group $\SL(n,\integer)$ itself is not left orderable, because it has elements of finite order. But there are subgroups of finite index that do not have any elements of finite order, and we will show that they, too, are not left orderable.
 \end{rem}

\begin{defn}
 Suppose $\prec$ is a left-invariant order relation on~$\Gamma$.
 For elements $a$~and~$b$ of~$\Gamma$, we say $a$ is
\defit{infinitely smaller} than~$b$ (denoted $a \ll
b$) if either $a^k \prec b$ for all $k \in \integer$ or
$a^k \prec b^{-1}$ for all $k \in\integer$. Notice that the
relation $\ll$ is transitive.
 \end{defn}

\begin{notation}
 The commutator $a ^{-1} b ^{-1} a b$ of
elements $a$~and~$b$ of a group~$\Gamma $ is denoted
$[a,b]$. It is straightforward to check that
$a$~commutes with~$b$ if and only if $[a,b] =
e$ (the identity element of~$\Gamma$).
 \end{notation}

\begin{lem} \label{CommutatorPowers}
 Let $a$~and~$b$ be elements of\/~$\Gamma$. If\/
$[a,b]$ commutes with both $a$ and~$b$, then\/
$[b^k,a^m] = [a,b]^{-km}$ for all $k,m \in
\integer$.
 \end{lem}
 
 \begin{proof}
 Exercise (or see \cite[Lems~2.2.2(i) and 2.2.4(iii), pp.~19 and~20]{GorensteinBook}).
 \end{proof}

\begin{lem}[(Ault, Rhemtulla)] \label{LOHeis}
 Suppose $a,b,z$ are non-identity elements of a
left-ordered group~$H$, with $[a,b] = z^k$ for some
nonzero $k \in \integer$, and $[a,z] = [b,z] = e$.
Then either $z \ll a$ or $z \ll b$.
 \end{lem}

\begin{proof}
 Assume, for simplicity, that $a,b,z \ge e$ and $k > 0$.
 (All other cases can be reduced to this one by replacing some or all of $a,b,z$ with their inverses, and/or interchanging $a$ with~$b$.)

We may
assume $z \not\ll a$, so $a \prec z^p$ for  some $p \in
\integer^+$. Similarly, we may assume $b \prec z^q$ for some $q
\in \integer^+$. Then (using the left-invariance of~$\prec$) we have
 $$ e \prec a^{-1} z^p,
 \qquad e \prec b^{-1} z^q,
 \qquad e \prec a,
 \qquad e \prec b,
 \qquad e \prec z .$$
 Hence, for all $m \in \integer^+$, we have
 \begin{align*}
 e &\prec  (b^{-1} z^q)^m(a^{-1} z^p)^m b^m a^m
  \\
 &= b^{-m} a^{-m} b^m a^m z^{(p+q)m}
 && \mbox{($z$ commutes with $a$ and~$b$)} \\
 &= [b^m, a^m] z^{(p+q)m}
 \\&= z^{-km^2} z^{(p+q)m}
 && \mbox{($[a,b] = z^k$ and see \pref{CommutatorPowers})} \\
 &= z^{-km^2 + (p+q)m} \\
 &= z^{\text{negative}}
 && \mbox{(if $m$ is sufficiently large)} \\
 &\prec e 
 .
 \end{align*}
 This is a contradiction.
 \end{proof}

\begin{proof}[\normalfont\textbf{Proof of \cref{SLnZnotLO}}]
 Suppose $\Gamma$ is left orderable. (This will lead to a
contradiction.)
 Because subgroups of left-orderable groups are left-orderable, we may assume $n = 3$; that is, $\Gamma $ has finite index in $\SL(3,\integer)$. Hence, there is
some $k \in \integer^+$ such that the six matrices
\begin{equation} \label{ElemUnipMatsInSL3}
 \begin{matrix}
 &a_1 = \begin{bmatrix} 1 & k & 0 \\ 0 & 1 & 0 \\ 0 & 0 & 1
\end{bmatrix},
 &a_2 = \begin{bmatrix} 1 & 0 & k \\ 0 & 1 & 0 \\ 0 & 0 & 1
\end{bmatrix},
 &a_3 = \begin{bmatrix} 1 & 0 & 0 \\ 0 & 1 & k \\ 0 & 0 & 1
\end{bmatrix}, \\
\ \\
 &a_4 = \begin{bmatrix} 1 & 0 & 0 \\ k & 1 & 0 \\ 0 & 0 & 1
\end{bmatrix},
 &a_5 = \begin{bmatrix} 1 & 0 & 0 \\ 0 & 1 & 0 \\ k & 0 & 1
\end{bmatrix},
 &a_6 = \begin{bmatrix} 1 & 0 & 0 \\ 0 & 1 & 0 \\ 0 & k & 1
\end{bmatrix} \hphantom{,}
\end{matrix}
 \end{equation}
 all belong to~$\Gamma$. A straightforward check shows that
$[a_i, a_{i+1}] = e$ and $[a_{i-1},a_{i+1}] =
(a_i)^{\pm k}$ for $i = 1,\dots,6$, with subscripts read
modulo~$6$. Thus, \cref{LOHeis} asserts
	\begin{equation} \label{SLnZnotLOPf-<<}
	\text{either $a_i \ll a_{i-1}$ or $a_i \ll a_{i+1}$.}
	\end{equation}
In particular, we must have either $a_1 \ll a_6$ or $a_1
\ll a_2$. Assume for definiteness that $a_1 \ll a_2$.
(The other case is very similar.) For each~$i$, \cref{LOHeis}
implies that if $a_{i-1} \ll a_i$, then $a_i \ll
a_{i+1}$. Since $a_1 \ll a_2$, we conclude by
induction that 
 $$ a_1 \ll a_2 \ll a_3 \ll a_4 \ll a_5 \ll
a_6 \ll a_1 .$$
 Thus $a_1 \ll a_1$, a contradiction.
 \end{proof}

For the reader acquainted with root systems,
let us give a more conceptual presentation of the proof of
\cref{SLnZnotLO}. 

\begin{proof}[{\textbf{Alternate version of the proof}}]
 The root system of $\SL(3,\integer)$ is of type~$A_2$, pictured in
\cref{RootsOfSL3Z}. For $i = 1,\ldots,6$, let $U_i$ be
the root space of~$\Gamma$ corresponding to the root~$\alpha_i$.
(Note that $U_i \iso \integer$ is cyclic. In the notation of the above
proof, the element~$a_i$ belongs to~$U_i$.) Define
 $$ \mbox{$\alpha_i \ll \alpha_j$ if there exists $u \in U_j$,
such that $U_i \prec u$} $$
 (that is, $v \prec u$, for all $v \in U_i$).

\begin{figure}
 \begin{center}
 \includegraphics{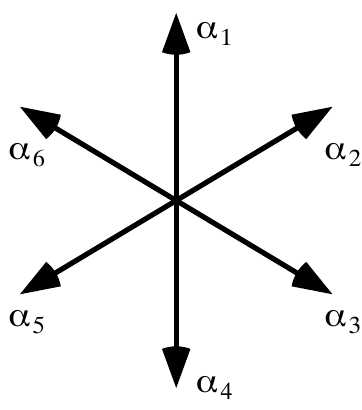}
 \caption{The root system of $\SL(3,\integer)$ (type $A_2$).}
 \label{RootsOfSL3Z}
 \end{center}
 \end{figure}

It is obvious, from \cref{RootsOfSL3Z}, that $\alpha_{i-1} +
\alpha_{i+1} = \alpha_i$, and that $\alpha_i + \alpha_{i+1}$ is not a
root, so 
 $$\mbox{$e \neq [U_{i-1}, U_{i+1}] \subset U_i$ 
 \qquad and \qquad
 $[U_i, U_{i+1}] = e$}
 .$$
 Hence, \cref{LOHeis} implies that either $\alpha_i \ll
\alpha_{i-1}$ or $\alpha_i \ll \alpha_{i+1}$. Arguing as in the above
proof, we conclude that
 $$ \alpha_1 \ll \alpha_2 \ll \alpha_3 \ll \alpha_4 \ll \alpha_5 \ll
\alpha_6 \ll \alpha_1 .$$
 Thus $\alpha_1 \ll \alpha_1$, a contradiction.
  \end{proof}

\begin{notes}
\Cref{SLnZnotLO} is due to D.\,Witte \cite{WitteCircle}.
(An exposition of the proof also appears in \cite[Thm.~7.2]{GhysCircleSurvey}.)
\Cref{LOHeis} is a special case of a theorem of J.\,C.\,Ault \cite{Ault-RONilp} and A.\,H.\,Rhemtulla \cite{Rhemtulla-ROGrps}:
if $\Lambda$ is any nontrivial, finitely generated, left-ordered,
nilpotent group, then there is some $a \in \Lambda$, such that
$[\Lambda,\Lambda] \prec a$.

An elementary proof that $\SL(n,\integer)$ has no nontrivial actions on the circle can be found in \cite{BridsonVogtmann-homos}, but the argument uses the existence of elements of finite order, so it does not apply to finite-index subgroups of $\SL(n,\integer)$.
\end{notes}

\section{The Reeb-Thurston Stability Theorem} \label{ReebThurstonSect}

The Ghys-Burger-Monod Theorem \pref{GhysBurgerMonodFinOrb} tells us that if $\Gamma$ is a lattice in $G = \SL(3,\real)$, then any action of $\Gamma$ on~$\circle$ has a finite orbit; in other words, some finite-index subgroup~$\Gamma'$ of~$\Gamma$ has a fixed point. In order to deduce \cref{LattSLnNoC1Act} from this, it suffices to show that if the action is by (orientation-preserving) $C^1$~diffeomorphisms, then $\Gamma'$ acts trivially. This triviality of~$\Gamma'$ is immediate from the following result:

\begin{prop}[({Reeb-Thurston Stability Theorem})]
\label{ReebThurston}
 Suppose
 \noprelistbreak
 \begin{itemize}
 \item $\Lambda$ is a finitely generated subgroup of
$\Diff^1_+(\circle)$,
 \item $\Lambda$ has a fixed point,
 and
 \item the abelianization $\Lambda/[\Lambda,\Lambda]$ is finite.
 \end{itemize}
 Then $\Lambda = \{e\}$ is trivial.
 \end{prop}

A differentiable action on~$\circle$ that has a fixed point can be
transformed into a differentiable action on the unit interval $[0,1]$
(cf.\ pf.\ of \ref{AOC<>LO}($\Rightarrow$)). Thus,
\cref{ReebThurston} can be reformulated as follows.

\begin{thmref}{ReebThurston}
 \begin{prop}[({Reeb-Thurston Stability Theorem})] \label{ReebThurstonI}
 Suppose
 \noprelistbreak
 \begin{itemize}
 \item $I = [0,1]$ is the unit interval,
 \item $\Lambda$ is a finitely generated subgroup of $\Diff^1_+(I)$,
 and
 \item the abelianization $\Lambda/[\Lambda,\Lambda]$ is finite.
 \end{itemize}
 Then $\Lambda = \{e\}$ is trivial.
 \end{prop}
 \end{thmref}

\begin{proof}[Proof in a special case]
 Define $\sigma \colon \Lambda \to \real^+$ by $\sigma(\lambda) =
\lambda'(0)$. From the Chain Rule, we see that $\sigma$ is a
(multiplicative) homomorphism. Because $\Lambda/[\Lambda,\Lambda]$ is
finite and $\real^+$ is abelian, this implies that $\sigma(\Lambda)$
is finite. However, $\real^+$ has no nontrivial, finite subgroups, so
this implies that $\sigma(\Lambda)$ is trivial; therefore
 $$ \mbox{$\lambda'(0) = 1$, for all $\lambda \in \Lambda$} .$$
 
 For simplicity, let us assume, henceforth, that each element
of~$\Lambda$ is real-analytic, rather than merely~$C^1$. (We
remark that there is no need to assume $\Lambda$ is finitely
generated in the real analytic case. The proof for $C^1$~actions is
similar, but requires some additional effort, and does use the hypothesis that $\Lambda$ is finitely generated.) Thus, each element~$\lambda$
of~$\Lambda$ can be expressed as a power series in a neighborhood
of~$0$:
 $$ \lambda(x) = x + a_{\lambda,2} x^2 + a_{\lambda,3} x^3 + \cdots
.$$
 (There is no constant term, because $\lambda(0) = 0$; the
coefficient of~$x$ is~$1$, because $\lambda'(0) = 1$.) 

Now suppose $\Lambda$ is nontrivial. (This will lead to a contradiction.) Then there exist $n$ and~$\lambda$, such that
	\begin{equation} \label{LambdaNontriv}
	a_{\lambda,n} \neq 0 
	. \end{equation}
We may assume $n$~is minimal; hence
	$$ \text{$\lambda(x) = x + a_{\lambda,n} x^n + a_{\lambda,n+1} x^{n+1} + \cdots$
	\ for all $ \lambda \in \Lambda$.} $$
Then, for
$\lambda,\gamma \in \Lambda$, we have
 \begin{align*}
 x + a_{\lambda \gamma,n} x^n + O(x^{n+1}) 
 &= (\lambda \gamma)(x) \\
 &= \lambda \bigl( \gamma(x) \bigr) \\
 &= \gamma(x) + a_{\lambda,n} \bigl( \gamma(x) \bigr)^n + O(x^{n+1}) \\
 &= \bigl( x + a_{\gamma,n} x^n + O(x^{n+1}) \bigr) 
 \\& \hskip 0.5in 
 + a_{\lambda,n}  \bigl( x + a_{\gamma,n} x^n + O(x^{n+1}) \bigr)^n
 \\& \hskip 0.9in
 + O(x^{n+1}) \\
 &= x + (a_{\gamma,n} + a_{\lambda,n}) x^n + O(x^{n+1}) 
 .
 \end{align*}
 Thus, the map $\tau \colon \Lambda \to \real$, defined by
$\tau(\lambda) = a_{\lambda,n}$, is an additive homomorphism. Because
$\Lambda/[\Lambda,\Lambda]$ is finite, but $\real$ has no nontrivial,
finite subgroups, this implies that $\tau(\Lambda)$ is trivial;
therefore $a_{\lambda,n} = 0$ for every $\lambda \in \Lambda$.
This contradicts \pref{LambdaNontriv}.
 \end{proof}

\begin{remgen} \label{ReebThurston-GenPf}
Suppose $\Lambda$ is nontrivial. Then there is no harm in assuming that $\Lambda$ acts nontrivially on every neighborhood of~$0$. 
Hence, letting $\Lambda_0$ be a finite generating set for~$\Lambda$, we may choose $\lambda_0 \in \Lambda_0$ and a sequence $x_n \to 0^+$, such that, for all~$n$, we have $\lambda_0(x_n) \neq x_n$ and
	$$ \text{$|\lambda(x_n) - x_n| \le |\lambda_0(x_n) - x_n|$ \ for all $\lambda \in \Lambda_0$}. $$
By passing to a subsequence of $\{x_n\}$, we may assume
	$$ \tau(\lambda) = \lim_{n \to \infty} \frac{\lambda(x_n) - x_n }{ \lambda_0(x_n) - x_n} $$
exists for all $\lambda \in \Lambda_0$. 
Because $\lambda'(0) = 1$ and $\lambda$~is $C^1$, it can be shown that $\tau(\lambda \gamma) = \tau(\lambda) + \tau(\gamma)$. Therefore $\tau\colon \Lambda \to \real$ is a (nontrivial) homomorphism. This is a contradiction.
\end{remgen}

\begin{rem}[(Zimmer)] \label{ReebThurston-ZimmerPf}
 Here is the outline of a nice proof of
\cref{ReebThurstonI}, under the additional assumption that
$\Lambda$~has Kazhdan's Property~$(T)$.
 \noprelistbreak
 \begin{enumerate}
 \item It suffices to show that the fixed points of~$\Gamma$ form a
dense subset of~$I$; thus, we may assume that $0$ and~$1$ are the
only points that are fixed by~$\Gamma$.
 \item We have $\gamma'(0) = 1$ for all $\gamma \in \Gamma$.
 \item Define a unitary representation of~$\Gamma$ on $L^2(I)$
by
 $$ f^\gamma(t) = f(\gamma t) \, |\gamma'(t)|^{1/2} .$$
 \item For any $\gamma \in \Gamma$ and any $\delta > 0$, if $f$~is
the characteristic function of a sufficiently small neighborhood
of~$0$, then $ \| f^\gamma - f \| < \delta \|f\|$. Hence, this unitary representation has almost-invariant vectors. 
 \item Because $\Lambda $ has Kazhdan's Property~$(T)$, this implies there are fixed vectors: there exists $f \in L^2(I) \smallsetminus \{0\}$, such that
$f^\gamma = f$ for all $\gamma \in \Gamma$.
 \item Every point in the essential support of~$f$ is a fixed point
of~$\Gamma$.
 \item This is a contradiction.
 \end{enumerate}
\end{rem}

\begin{notes}
\Cref{ReebThurston} was proved by W.\,Thurston
\cite{ThurstonStab} in a more general form that also applies to
actions on manifolds of higher dimension. (It generalizes a 
theorem of G.\,Reeb.) 

See \cite{ReebSchweitzer} or \cite{Schachermayer} for details of the proof sketched in \cref{ReebThurston-GenPf}.

The proof
outlined in \cref{ReebThurston-ZimmerPf} is due to R.\,J.\,Zimmer; details appear in
\cite[\S5, pp.~108--109]{WitteZimmer-ActOnCircle}.
\end{notes}

\section{Smooth actions of Kazhdan groups on the circle} \label{NavasPfSect}

In this section, we prove \cref{Navas-AOC}. First, let us recall one of the many equivalent definitions of Kazhdan's property $(T)$.

\begin{notation}
For any real Hilbert space~$\hilbert$, we use $\Isom(\hilbert)$ to denote the isometry group of~$\hilbert$. (We remark that each isometry of~$\hilbert$ is the composition of a translation with a norm-preserving linear transformation.)
\end{notation}

\begin{defn} \label{KazhdanDefn}
We say that a discrete group~$\Gamma$ has \defit{Kazhdan's property~$(T)$} if, for every homomorphism $\rho \colon \Gamma \to \Isom(\hilbert)$, where $\hilbert$ is a real Hilbert space, there exists $v \in \hilbert$, such that $\rho(g) v = v$, for every $g \in \Gamma$.
 \end{defn}
 
 In short, to say that $\Gamma$ has Kazhdan's property~$(T)$ means that every isometric action of~$\Gamma$ on any Hilbert space has a fixed point.
 The importance of this notion for our purposes stems from the following result, whose proof we omit.

\begin{thm}[(Kazhdan)]
If $n \ge 3$, then every lattice in $\SL(n,\real)$ has Kazhdan's property~$(T)$.
\end{thm}

\begin{exer}
Show that if $\Gamma$ is an infinite group with Kazhdan's property~$(T)$, then $\Gamma$ is not abelian.
\hint{Every group with Kazhdan's property~$(T)$ is finitely generated, and every finitely generated abelian group is either finite or has a quotient isomorphic to~$\integer$.}
\end{exer}

We now turn to the proof of \cref{Navas-AOC}. To simplify notation, 
we may think of $\circle$ as $[-\pi/2,\pi/2]$.
 In particular, for any diffeomorphism of~$\circle$ and any $x \in \circle$, the
derivative $g'(x)$ is a well-defined real number.

\begin{defn} \ 
\noprelistbreak
 \begin{itemize}
 \item Let $\Func(\circle \times \circle)$ be the vector space of
measurable functions on $\circle \times \circle$ (with two functions
being identified if they are equal almost everywhere).
 \item Define an action of $\Diff^2(\circle)$ on $\Func(\circle
\times \circle)$ by
 $$F^g(x,y) = F \bigl(g(x), g(y) \bigr)\, |g'(x)|^{1/2}
\,|g'(y)|^{1/2} $$
  for $F \in \Func(\circle \times \circle)$ and  $g \in
\Diff^2(\circle)$.
 \item Let $\|F\|_2 = \left( \int_{\circle} \int_{\circle} F(x,y)^2 \,
dx \, dy \right)^{1/2}$ be the $L^2$-norm of~$F$; note that $\|F\|_2 =
\infty$ if $F \notin L^2(\circle \times \circle)$.
 \end{itemize}
 Note that
$$ 
 \mbox{$F^{gh} = (F^g)^h$
 \quad
 and
 \quad
 $\|F^g\|_2 = \|F\|_2$
 \quad for $F \in \Func(\circle \times \circle)$ and  $g, h \in
\Diff^2(\circle)$.}
$$ 
\end{defn}

\begin{notation} \ 
\noprelistbreak
\begin{itemize}
\item Choose a positive function~$f$ on $\circle$, such that 
 \noprelistbreak
 \begin{itemize}
 \item $f$~has a $1/x$~singularity at the point~$0$ of~$\circle$, and
 \item $f$ is~$C^\infty$ everywhere else;
 \end{itemize}
 that is, identifying $\circle$ with $[-\pi/2,\pi/2]$, we have
 \begin{equation} \label{Navas-pole}  
 f(x) - \frac{1}{|x|} \in C^\infty(\circle) 
 .
 \end{equation}
 For example, one may take $f(x) = | \cot x  \, |$.

\item Now define 
 $$ \mbox{$\Phi(x,y) = f(x-y)$ on $\circle \times \circle$.} $$
\end{itemize}
 Because of the $1/x$ singularity of~$f$, it is easy to see that
$\Phi \notin L^2(\circle \times \circle)$.
\end{notation}

For any $g \in \Diff^2(\circle)$, the following calculation shows
that the singularity of
$\Phi^g$ cancels the singularity of~$\Phi$.

\begin{lem} \label{NavasBddDiff}
The
difference $\Phi^g - \Phi$ is a bounded function on $\circle \times
\circle$, for any $g \in \Diff^2(\circle)$.
 \end{lem}

\begin{proof}
 From \pref{Navas-pole}, we see that there is no harm in
working with the function $\Phi_0(x,y) = 1/|x-y|$, instead of~$\Phi$.
Also, in order to reduce the number of absolute-value signs, let us assume $g' \ge 0$ everywhere.
 \begin{align*}
& |\Phi^g(x,y) - \Phi(x,y)| 
 \\&\qquad \approx |\Phi_0^g(x,y) - \Phi_0(x,y)| \\
 &\qquad = \left| g'(x)^{1/2} g'(y)^{1/2} \, \Phi_0 \bigl( g(x), g(y) \bigr)
 - \Phi_0(x,y) \right| \\
 &\qquad = \left| \frac{g'(x)^{1/2} g'(y)^{1/2} }{| g(x) - g(y) |}
 - \frac{1}{|x-y|} \right| \\
 &\qquad =  \left| \frac{g'(x)^{1/2} g'(y)^{1/2} }{ \bigl| g'(t) (x-y) \bigr|}
 - \frac{1}{|x-y|} \right| 
 && \begin{pmatrix} \text{$\exists t \in (x,y)$, by} \\ \text{Mean Value Thm.} \end{pmatrix} \\
 &\qquad = \frac{ | g'(x)^{1/2} g'(y)^{1/2} - g'(t)| }{ g'(t) \, |x-y|} \\
 &\qquad = \frac{ | g'(x) g'(y) - g'(t)^2| }
 { \bigl( g'(x)^{1/2} g'(y)^{1/2} + g'(t) \bigr) g'(t)|x-y|
  }
 && \begin{pmatrix} \text{multiply by conjugate} \\ \text{of numerator} \end{pmatrix} \\
 &\qquad = O \left( \frac{ | g'(x) g'(y) - g'(t)^2| }{ |x-y|}\right) 
 && \begin{pmatrix} \text{$g$ diffeomorphism} \\ \text{$\Rightarrow$ $g'$ is never~$0$} \end{pmatrix} 
 .
 \end{align*}
 Now, we wish to show that the numerator is bounded by a constant
multiple of the denominator.
 \begin{align*}
 &| g'(x) g'(y) - g'(t)^2| 
 \\&\qquad \le g'(x) | g'(y) - g'(t)| + g'(t)|g'(x) - g'(t)| 
  && \mbox{(Triangle inequality)} \\
 &\qquad = g'(x)\, |g''(u)| \,|y-t|   + g'(t)\, |g''(v)| \, |x-t| 
 && \begin{pmatrix} \text{Mean Value Thm.} \\ \text{applied to~$g'$} \end{pmatrix} \\
 &\qquad \le g'(x)\, |g''(u)| \,|x-y| + g'(t)\, |g''(v)| \, |x-y|
 && \mbox{(because $t \in (x,y)$)} \\
 &\qquad = O \bigl( |x-y| \bigr)
 && \begin{pmatrix} \text{$g',g''$ continuous,} \\ \text{so bounded} \end{pmatrix}
 .
 \end{align*}
 \end{proof}
 
 We will also use the following classical fact:
 
 \begin{lem}[(H\"older, 1901)] \label{Nonab->FP}
 Every nonabelian group of homeomorphisms of~$\circle$ contains at least one nonidentity element that has a fixed point.
 \end{lem}

\begin{proof}[\normalfont\textbf{Proof of \cref{Navas-AOC}}]
Suppose $\Gamma$ is a discrete group with Kazhdan's property~$(T)$, and we have a faithful $C^2$ action of~$\Gamma$ on the circle~$\circle$. From \cref{NavasBddDiff}, we see that
 \begin{equation} \label{Navas-Hinv}
 \mbox{$\bigl(\Phi + L^2(\circle \times \circle) \bigr)^g = \Phi +
L^2(\circle \times \circle)$ for all $g \in \Diff^2(\circle)$.}
 \end{equation}
 Thus, $\Gamma$~acts (by isometries) on the affine Hilbert space
$\Phi + L^2(\circle \times \circle)$.
Because $\Gamma$ has Kazhdan's Property, 
we conclude that $\Gamma$ has a fixed point~$F$ in 
 $\Phi + L^2(\circle \times \circle)$:
 $$ \mbox{$F^g = F$ for all $g \in \Gamma$} .$$
 Because $\Phi \notin L^2(\circle \times \circle)$, and $F - \Phi \in
L^2(\circle \times \circle)$, it is obvious that $F \notin
L^2(\circle \times \circle)$.

Now, define a measure~$\mu$ on~$\circle \times \circle$ by $\mu = F^2
\, dx \, dy$. Then
 \noprelistbreak
 \begin{enumerate}
 \item $\mu$ is a $\Gamma$-invariant measure on $\circle \times
\circle$ (because $F$~is $\Gamma$-invariant);
 and
 \item if $R = (a_1,b_1) \times (a_2,b_2)$ is a rectangle in $\circle
\times \circle$, then, from the $1/x$ singularity of~$\Phi$ along the diagonal, we see that
 $$\mu(R) =
 \begin{cases}
 \infty & \mbox{if $R$ intersects the diagonal;} \\
 \mbox{finite} & \mbox{if $R$ is away from the diagonal}
 .
 \end{cases}$$
 For example, in \cref{NavasRectabcFig}, the rectangle $(a,b)
\times (b,c)$ has infinite measure,
because it touches the diagonal at $(b,b)$, but the shaded rectangle
$R_k$ has finite measure, because it does
not touch the diagonal.
 \end{enumerate}

\begin{figure}
\begin{center}
\includegraphics[scale=0.500]{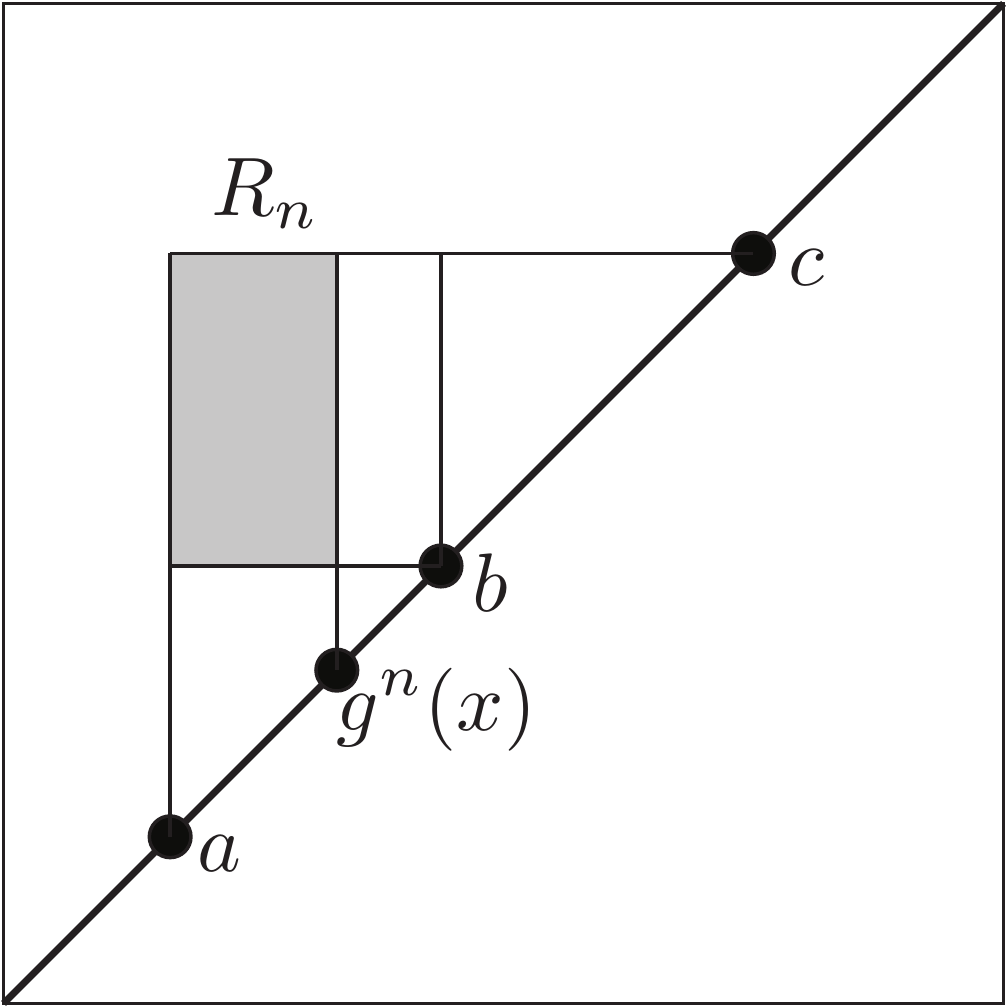}
 \caption{A rectangle $R_n$ (shaded) in $\circle \times \circle$.}
 \label{NavasRectabcFig}
 \end{center}
 \end{figure}

Because $\Gamma$ is not abelian, \cref{Nonab->FP} tells us     we may choose $g \in \Gamma$, such
that $g$~has a fixed point.  Then, by
passing to a triple cover of~$\circle$, we obtain an action of (a
finite extension of)~$\Gamma$ in which $g$ has at least~$3$ fixed
points:
 $$ \mbox{$g(a) = a$, \qquad $g(b) = b$, \qquad $g(c) = c$.}$$
 By perhaps replacing $a$~and~$b$ with different fixed points, we may
assume, for $x \in (a,b)$, that
 $$\lim g^k(x) =
 \begin{cases}
 b & \mbox{as $k \to \infty$} \\
 a & \hbox{as $k \to -\infty$}
 .
 \end{cases}$$

For each $k \in \integer$, define 
 $$R_k = \bigl( a, g^k(x) \bigr) \times (b,c) $$
 \csee{NavasRectabcFig}.
 Then
 \begin{align*}
 g(R_k) &= R_{k+1}, \\
 \intertext{so}
  \mu(R_k) &= \mu(R_{k+1})
 .
 \end{align*}

Therefore
 \begin{align*}
 0 &= \sum_{k=-\infty}^\infty \bigl( \mu(R_{k+1}) - \mu(R_k) \bigr) \\
 &= \mu\left( \bigcup_{k=-\infty}^\infty (R_{k+1} \smallsetminus R_k)
\right) \\
 &= \mu\bigl( (a,b) \times (b,c) \bigr)
 .
 \end{align*}

However, the rectangle $(a,b) \times (b,c)$ touches the diagonal at
the point $(b,b)$, so it has infinite measure. This is a contradiction.
\end{proof}

\begin{rem} \label{Navas-AOC3/2Rem}
 In the proof of \cref{Navas-AOC}, the assumption that the
elements of $\Gamma$ are~$C^2$ was used only to establish
\pref{Navas-Hinv}. For this, it is not necessary to show that
$\Phi^g - \Phi$ is bounded (as in \pref{NavasBddDiff}), but only that
$\Phi^g - \Phi \in L^2(\circle \times \circle)$. The calculations in the proof of \cref{NavasBddDiff} show that this holds under the
weaker hypothesis that $g \in C^{3/2 + \epsilon}$, for any $\epsilon > 0$.
In fact, A.\,Navas observed that, by using recent results that groups with Kazhdan's property~$(T)$ also have fixed points in certain $L^p$ spaces with $p \neq 2$, it can be shown that  $C^2$ can be replaced with $C^{3/2}$ in the statement of \cref{Navas-AOC}.
 \end{rem}

This leads to the following well-known question:

\begin{problem}
 Can an infinite (discrete) group with Kazhdan's property $(T)$ have a faithful $C^0$ action on~$\circle$?
  \end{problem}

 \begin{notes}
 \Cref{Navas-AOC} was first proved in \cite{NavasCircle},
 but the idea to use \cref{NavasBddDiff} came from earlier work of A.~Pressley
and G.~Segal \cite{PressleySegal} and A.~Reznikov
\cite[Chap.~2]{Reznikov}. Proofs also appear in \cite[\S2.9]{BekkaEtAlT} and \cite[\S5.2]{Navas-GrpsDiffeos}.

 The condition in \Cref{KazhdanDefn} was introduced by J.--P.\,Serre, and it is not at all obvious that it is equivalent to the original definition of property~$(T)$ that was given by D.\,Kazhdan.
  For a discussion of this, and much more, the standard reference on Kazhdan's property~$(T)$ is \cite{BekkaEtAlT}.

 See \cite[Thm.~6.10]{GhysCircleSurvey} or \cite[\S2.2.4]{Navas-GrpsDiffeos} for a proof of \cref{Nonab->FP}.
 
Implications of Kazhdan's property~$(T)$ for fixed points of actions on $L^p$ spaces (and other Banach spaces) are discussed in \cite{BaderEtAlT}.
 \end{notes}

\section{Ghys's proof that actions have a finite orbit} \label{GhysPfSect}

In this section, we present Ghys's proof of \cref{GhysBurgerMonodFinOrb}, modulo some facts that will be proved in \cref{GhysMissingSect}.
To get started, let us show that it suffices to find a $\Gamma$-invariant
measure on the circle.

\begin{defn}
A measure~$\mu$ on a measure space~$X$ is a \defit{probability measure} if $\mu(X) = 1$.
\end{defn}

\begin{lem} \label{GhysInvtMeasSuff}
 If 
 \noprelistbreak
 \begin{itemize}
 \item $\Gamma$ is a discrete group, such that the abelianization\/ $\Gamma/[\Gamma,\Gamma]$ is finite,
 \item $\Gamma$ acts on~$\circle$ by orientation-preserving homeomorphisms,
 and
 \item there is a\/ $\Gamma$-invariant probability measure~$\mu$
on~$\circle$,
 \end{itemize}
 then\/ $\Gamma$ has at least one finite orbit on~$\circle$.
 \end{lem} 

\begin{proof}
 We consider two cases.

\setcounter{case}{0}

\begin{case}
 Assume $\mu$ has at least one atom~$p$.
 \end{case}
 Let 
 $$ \Gamma p = \{\, g p \mid g \in \Gamma \,\} $$
 be the orbit of~$p$. Because $\mu$ is $\Gamma$-invariant, we
have
 $\mu \bigl( g p \bigr) = \mu(p)$,
 for every $g \in \Gamma$. Therefore
 $$\mu(\Gamma p)
 = \# (\Gamma p) \, \mu(p) .$$
 Since $\mu$ is a probability measure, we know that $\mu(\Gamma
p) < \infty$, so we conclude that $\Gamma p$~is a finite
set. That is, the orbit of~$p$ is finite.

\begin{case} \label{GhysInvtMeasSuffPf-noatoms}
 Assume $\mu$ has no atoms.
 \end{case}
 To simplify the proof, let us assume that the support of~$\mu$ is
all of~$\circle$. (In other words, every
nonempty open subset of~$\circle$ has positive measure.) For $x,y \in
\circle$, define
 \begin{equation} \label{GhysFPSuff-DefMetric}
  d(x,y) = \mu \bigl( [x,y] \bigr) ,
 \end{equation}
 where $[x,y]$ is a path from~$x$ to~$y$, and, for a given $x$
and~$y$, we choose the path $[x,y]$ to minimize $\mu \bigl( [x,y]
\bigr)$.
 It is easy to see that $d$~is a metric on~$\circle$.
  Up to isometry, there is a unique metric
on~$\circle$, so we may assume that $d$~is
the usual arc-length metric.

Because $\mu$ is $\Gamma$-invariant, we know that $d$~is
$\Gamma$-invariant, so $\Gamma$ acts by rotations of the circle.
There is no harm in assuming that the action is faithful, so we conclude that $\Gamma$ is abelian. But every abelian quotient
of~$\Gamma$ is finite, so we conclude that
$\Gamma$ is finite. Hence, every orbit is finite.
 \end{proof}
 
 \begin{exer}
 Complete the proof of \cref{GhysInvtMeasSuff}, by eliminating the assumption that the support of~$\mu$ is all of~$\circle$ in \cref{GhysInvtMeasSuffPf-noatoms}.
 \hint{Modify the above proof to show that every orbit in the support of~$\mu$ is finite.}
 \end{exer}

To simplify the notation, we will assume henceforth that $n = 3$.

\begin{notation} \label{SL3Notation}
Let
\noprelistbreak
\begin{itemize}
\item $G = \SL(3,\real)$,
\item $\Gamma$ be a lattice in~$G$,
and
 \item $P = \begin{bmatrix}
 * & * & * \\
  & * & * \\
  &  & *
 \end{bmatrix}
 \subset G$
 \end{itemize}
 We remark that $P$ is a minimal parabolic subgroup of~$G$. 
 \end{notation}
 
  Ghys's proof is based on the following key fact, which will be established in \cref{BdryMapSubsec} by using the fact that the group~$P$ is ``amenable\..''
 
\begin{notation}
$\Prob(X)$ denotes the set of all Radon probability measures on~$X$ (where $X$ is any compact, Hausdorff space).
This is a closed, convex subset of the unit ball in $C(X)^*$, with the weak$^*$ topology, so $\Prob(X)$ has a natural topology that makes it a compact Hausdorff space.
\end{notation}

 \begin{prop}[(Furstenberg)] \label{BdryMap}
 If $\Gamma$ acts continuously on any compact metric space~$X$, then there is a $\Gamma$-equivariant, Borel measurable map $\psi \colon G/P \to \Prob(X)$.
 \end{prop}

 If $\psi$ were \emph{in}variant, rather than \emph{equi}variant, the following fundamental theorem would immediately imply that $\psi$ is constant (a.e.). (This theorem will be proved in \cref{MooreErgPfSect}.)
 
 \begin{thm}[(Moore Ergodicity Theorem)] \label{MooreErgThm}
 If $H$ is any closed, noncompact subgroup of~$G$, then the action of\/~$\Gamma$ on $G/H$ is \defit{ergodic}: by definition, this means that every\/ $\Gamma$-invariant, measurable function on $G/H$ is constant\/ {\upshape(}a.e.{\upshape)}.
 \end{thm}
 
 The homogeneous space $G/P$ plays a major role in Ghys's proof, because of its appearance in \cref{BdryMap}. We will now present a geometric interpretation of this space that  is very helpful.

Recall that a \defit{flag} in~$\real^3$ consists
of a pair $(\ell, \pi)$, where 
 \noprelistbreak
 \begin{itemize}
 \item $\ell$~is a $1$-dimensional vector subspace of~$\real^3$ (a
line), and 
 \item $\pi$ is a $2$-dimensional vector subspace of~$\real^3$ (a
plane) that contains~$\ell$.
 \end{itemize}
 The group $G = \SL(3,\real)$ acts transitively on the set of flags,
and $P$~is the stabilizer of the standard flag
 $$ \mbox{$F_0 = (\ell_0, \pi_0)$,
  where $\ell_0 = (*,0,0)$
 and $\pi_0 = (*,*,0)$}. $$
 Therefore:
 
 \begin{prop}
 $G/P$ can be identified with the space~$\flags$ of all flags, by
identifying $g P$ with the flag $g F_0$. 
\end{prop}

\begin{proof}[{\normalfont\textbf{Proof of \cref{GhysBurgerMonodFinOrb}}}]
Suppose $\Gamma$ acts by homeomorphisms on the circle~$\circle$.
 From \cref{BdryMap}, we know
there is a $\Gamma$-equivariant, measurable map 
 $$ \psi \colon G/P \to \Prob(\circle) . $$
 It will suffice to show that $\psi$ is constant (a.e.), for
then:
 \noprelistbreak
 \begin{itemize}
 \item the essential range of~$\psi$ consists of a single point $\mu
\in \Prob(\circle)$,
 and
 \item the essential range of~$\psi$ is $\Gamma$-invariant, because
$\psi$~is $\Gamma$-equivariant.
 \end{itemize} 
 So $\mu$ is $\Gamma$-invariant, and then \cref{GhysInvtMeasSuff} implies $\Gamma$ has a finite orbit on~$\circle$.

 There are two basic cases to consider: either the measure $\psi(x)$
consists entirely of atoms, or $\psi(x)$ has no atoms.
(Recall that an \defit{atom} of a measure~$\mu$ is a point~$p$, such that $\mu\bigl( \{p\} \bigr) \neq 0$.) It is also
possible that $\psi(x)$ consists partly of atoms, and partly of
nonatoms, but \cref{EquiPsiAtoms} below tells us that we need only consider the two extreme cases.

\setcounter{case}{0}

\begin{case} \label{GhysThmPf-noatoms}
 Assume $\psi(x)$ has no atoms, for a.e.~$x \in G/P$.
 \end{case}
 \noprelistbreak
 \begin{itemize}
 \item Let
 $ \Prob_0(\circle) = \{\, \mu \in \Prob(\circle) \mid
 \mbox{$\mu$ has no atoms} \,\} $.
 \item By assumption, we know $\psi \colon G/P \to \Prob_0(\circle)$,
so we may define
 $$ \text{$\psi^2 \colon (G/P)^2 \to \bigl( \Prob_0(\circle)  \bigr)^2$
 \quad  by \quad
 $\psi^2(x,y) = \bigl( \psi(x), \psi(y) \bigr)$.} $$
 Then $\psi^2$,
like~$\psi$, is measurable and $\Gamma$-equivariant.
 \item Define $d \colon \bigl( \Prob_0(\circle)  \bigr)^2 \to \real$
by
 $$ d(\mu_1,\mu_2) = \sup_J |\mu_1(J) - \mu_2(J)| ,$$
 where $J$~ranges over all intervals (that is, over all connected
subsets of~$\circle$). 
Since $\Gamma$ acts by homeomorphisms, and any homeomorphism of~$\circle$ maps intervals to intervals, it is easy to see that $d$~is $\Gamma$-invariant.
 \end{itemize}

 We claim that
 \begin{equation} \label{GhysPf-dCont}
 \mbox{$d$~is continuous,}
 \end{equation}
 with respect to the usual weak$^*$ topology that $\Prob_0(\circle)$
inherits from being a subset of $\Prob(\circle)$.
To see this, we note that if $\mu \in \Prob_0(\circle)$ and $\epsilon > 0$, then, because $\mu$~has no atoms, we may partition~$\circle$ into finitely many intervals $J_1,\ldots,J_k$, such that $\mu(J_i) < \epsilon$ for all~$i$.
By approximating the characteristic functions of these intervals from above and below, we may construct continuous functions $f_1,\ldots,f_{2k}$ and some $\delta > 0$, such that, for $\nu \in \Prob(\circle)$,
	$$ \text{if $|\nu(f_i) - \mu(f_i)| < \delta$ for every~$i$, then $|\nu(J_i) - \mu(J_i)| < \frac{\epsilon}{n}$ for every~$i$.} $$
Then, for any interval~$J$ in~$\circle$, we have $|\nu(J) - \mu(J)| < 2 \epsilon$, so $d(\nu,\mu) < 2 \epsilon$. This completes the proof of~\pref{GhysPf-dCont}.

Because $\psi^2$ is measurable, and $d$~is continuous, we know that the composition $d
\circ \psi^2$ is measurable.
In addition, because $\psi^2$~is $\Gamma$-equivariant and
$d$~is $\Gamma$-invariant, we also know that $d
\circ \psi^2$ is $\Gamma$-invariant.  (\emph{Note:} we are saying
\emph{in}variant, not just \emph{equi}variant.) From the Moore
Ergodicity Theorem \pref{MooreErgThm}, we know that $\Gamma$ is ergodic on $(G/P)^2$
\see{ErgOn(G/P)2}, so we conclude that $d \circ \psi^2$ is
 constant (a.e.); say 
 	$$ \text{$d \bigl( \psi^2(x,y) \bigr) = c$, \quad for
a.e.~$x,y \in G/P$.} $$
We wish to show that $c = 0$, for then it is clear that $\psi(x) =
\psi(y)$ for a.e.~$x,y \in G/P$, so $\psi$~is  constant (a.e.).

It is obvious that $d \bigl( \psi^2(x,x) \bigr) = 0$ for every $x \in
G/P$. If $\psi$ were constant everywhere, rather than merely
\emph{almost} everywhere, then it would follow immediately that $c
= 0$. Unfortunately, the diagonal $\{(x,x)\}$ is a set of measure~$0$
in $(G/P)^2$, so a little bit of additional argument is required.

By Lusin's Theorem, $\psi$~is continuous on a set~$K$ of positive measure.
Since the composition of continuous functions is continuous, we conclude that $d \circ \psi^2$ is continuous on $K \times K$. So, by continuity, it is constant on all of $K \times K$, not merely almost all. Since $d \bigl( \psi^2(x,x) \bigr) = 0$, this implies $d \bigl( \psi^2(x,y) \bigr) = 0$ for all $x,y \in K$. Since $K$~is a set of positive measure, this implies $c = 0$, as desired.

\begin{case} \label{GhysThmPf-atoms}
 Assume $\psi(x)$ consists entirely of atoms, for a.e.~$x \in G/P$.
 \end{case}
 To simplify the notation, without losing the main ideas, let us
assume that $\psi(x)$ consists of a single atom, for every $x \in G/P$. Thus, we may think of~$\psi$ as a
$\Gamma$-equivariant, measurable map
 $$ \psi \colon G/P \to \circle .$$
 Surprisingly, even with the simplifying assumption, the argument here
seems to be more difficult than in \cref{GhysThmPf-noatoms}.
The idea is to obtain a contradiction from the $\Gamma$-equivariance
of~$\psi$, by contrasting two fundamental observations:
 \noprelistbreak
 \begin{itemize}
 \item $\Homeo_+(\circle)$ is \emph{not} triply transitive
on~$\circle$: if $x$, $y$, and~$z$ are distinct, then no
orientation-preserving homeomorphism of~$\circle$ can map the triple
$(x,y,z)$ to $(y,x,z)$ --- they have opposite orientations under the
circular order on~$\circle$.
 \item The action of $\GL(2,\real)$ on the projective line $\real P^1
= \real \cup \{\infty\}$ by linear-fractional transformations
 $$ g(x) = \frac{ax+b}{cx+d} 
 \mbox{ \qquad if $g = \begin{bmatrix}
 a & b \\
 c & d
 \end{bmatrix}$} $$
 \emph{is} triply transitive: if $(x_1, y_1, z_1)$ and $(x_2, y_2,
z_2)$ are two ordered triples of distinct elements of $\real \cup
\{\infty\}$, then there is some $g \in  \GL(2,\real)$ with $g(x_1) =
x_2$, $g(y_1) = y_2$, and $g(z_1) = z_2$.
 \end{itemize}

To illustrate, let us give an easy proof that is not quite correct;
the actual proof is a modified version of this. Define
 $$ \psi^3 \colon (G/P)^3 \to (\circle)^3
 \mbox{ \qquad by \qquad $\psi^3(x,y,z) = \bigl( \psi(x), \psi(y), \psi(z)
\bigr)$} .$$
 Then $\psi^3$ is $\Gamma$-equivariant, so
 $$ X^+ = \{\, (x,y,z) \in (G/P)^3 \mid
 \mbox{$\bigl( \psi(x), \psi(y), \psi(z) \bigr)$ is positively
oriented} \,\} $$
 is a $\Gamma$-invariant, measurable subset of $(G/P)^3$. Let us
assume that $\Gamma$ is ergodic on $(G/P)^3$. (Unfortunately, this assumption is
false, so it is where the proof breaks down.)
Then $X^+$ must be (almost) all of $(G/P)^3$; thus, $\bigl(
\psi(x), \psi(y), \psi(z) \bigr)$ is positively oriented, for
(almost) every $x,y,z \in G/P$. But this is nonsense: either $\bigl(
\psi(x), \psi(y), \psi(z) \bigr)$ or $\bigl( \psi(y), \psi(x),
\psi(z) \bigr)$ is negatively oriented, so there are many negatively
oriented triples.

To salvage the above faulty proof, we replace $(G/P)^3$ with a
subset~$X$, on which $\Gamma$ does act ergodically. Let
 \noprelistbreak
 \begin{itemize}
 \item $Q = 
 \begin{bmatrix}
 * & * & * \\
 * & * & * \\
 0 & 0 & * 
 \end{bmatrix}
 \subset G$,
 and
 \item $X = 
 \bigset{
 (x_1, x_2, x_3)
 \in (G/P)^3
 }{
 \begin{matrix}
 x_1 Q = x_2 Q = x_3 Q, \\
 \mbox{$x_1, x_2, x_3$ distinct}
 \end{matrix}
 }$.
 \end{itemize}
 Note that $X$ is a submanifold of $(G/P)^3$. (If this is not
obvious, it follows from the fact, proven below, that $X$~is a single
$G$-orbit in $(G/P)^3$.)

 Assume, for the moment, that $\Gamma$ is ergodic on~$X$ (with
respect to any (hence, every) Lebesgue measure on the manifold~$X$).
Then the above proof, with $X$~in the place of $(G/P)^3$, implies that
 $$ X^+ = \{\, (x,y,z) \in X \mid
 \mbox{$\bigl( \psi(x), \psi(y), \psi(z) \bigr)$ is positively
oriented} \,\} $$
 is a set of measure~$0$. Then it is not difficult to see that $\psi$ is
 constant on~$X$ (a.e.). Hence $\psi$ is
right $Q$-invariant (a.e.): for each $q \in Q$, we have $\psi(xqP) =
\psi(xP)$ for a.e.\ $x \in G/P$.

By a similar argument, we see that $\psi$ is  right
$Q'$-invariant (a.e.), where
 $$ Q' = 
 \begin{bmatrix}
 * & * & * \\
 0 & * & * \\
 0 & * & * 
 \end{bmatrix}
 \subset G .$$
 Because $Q$~and~$Q'$, taken together, generate all of~$G$, it then follows that $\psi$ is  right
$G$-invariant  (a.e.). Hence, $\psi$ is 
constant (a.e.), as desired.

All that remains is to show that $\Gamma$ is ergodic on~$X$. By the
Moore Ergodicity Theorem \pref{MooreErgThm}, we need only show
that 
 \noprelistbreak
 \begin{enumerate}
 \item \label{GhysThmPf-GtransX}
 $G$ is transitive on~$X$, and
 \item \label{GhysThmPf-Stab}
 the stabilizer of some element of~$X$ is not compact.
 \end{enumerate}
 These facts are perhaps easiest to establish from a geometric perspective.

The (parabolic) subgroup~$Q$ is the stabilizer of the plane $\pi_0 = (*,*,0)$.
Hence, $Q F_0$ is the set of all flags $(\ell,\pi)$ with $\pi = \pi_0$.
Therefore, under the identification of $G/P$ with~$\flags$, we have
 \begin{equation} \label{GhysPf-X=flags}
 X = \bigset{
 \bigl( (\ell_1, \pi_1), (\ell_2, \pi_2), (\ell_3, \pi_3) \bigr)
 \in \flags^3
 }{
 \begin{matrix}
 \pi_1 = \pi_2 = \pi_3, \\
 \mbox{$\ell_1, \ell_2, \ell_3$ distinct}
 \end{matrix}
 }
. \end{equation}

\pref{GhysThmPf-GtransX}
 Let us show that $G$ is transitive on~$X$. Given 
 $$ \bigl( (\ell_1, \pi), (\ell_2, \pi), (\ell_3, \pi) \bigr) ,
 \bigl( (\ell_1', \pi'), (\ell_2', \pi'), (\ell_3, \pi') \bigr) \in X
,$$
 it suffices to show that there exists $g \in G$, such that 
 \begin{equation} \label{GhysThmPf-g(a,b,c)=a',b',c')}
 g \bigl( (\ell_1, \pi), (\ell_2, \pi), (\ell_3, \pi) \bigr) 
 = \bigl( (\ell_1', \pi'), (\ell_2', \pi'), (\ell_3, \pi') \bigr) .
 \end{equation}
 Because $G$ is transitive on the set of $2$-dimensional subspaces,
we may assume 
 $$\pi = \pi' = \pi_0 = \real^2 .$$
 Then, because $\GL(2,\real)$ is triply transitive on~$\real P^1$,
there exists $T \in \GL(2,\real)$, such that 
 $$T(\ell_1,\ell_2,\ell_3) = (\ell_1',\ell_2',\ell_3') .$$
 Letting
 $$ g = 
 \begin{bmatrix}
\vbox to 0pt{\vskip-6pt\hbox to 0pt{\hskip4pt\Huge $T$\hss}\vss} & & 0 \\
 && 0 \\
 0&0 & \frac{1}{\det T}
 \end{bmatrix} $$
 yields \pref{GhysThmPf-g(a,b,c)=a',b',c')}.

 \pref{GhysThmPf-Stab} 
 Let us show that the stabilizer of some element of~$X$ is not
compact. 
 Let 
 $$ A = 
 \bigset{
 \begin{bmatrix}
 a & 0 & 0 \\
 0 & a & 0 \\
 0 & 0 & 1/a^2
 \end{bmatrix}
 }{
 a \in \real^\times
 }
 .$$
 Then $A$ is a closed, noncompact subgroup of~$G$.  Furthermore,
every element of~$A$ acts as a scalar on~$\pi_0$, so every element
of~$A$ fixes every $1$-dimensional vector subspace of~$\pi_0$. Thus,
if
 $$ \bigl( (\ell_1, \pi), (\ell_2, \pi), (\ell_3, \pi) \bigr) \in X
,$$
 with $\pi = \pi_0$, then $A$ is contained in the stabilizer of this
element of~$X$, so the stabilizer is not compact.
 \end{proof}

\begin{notes}
Ghys's proof first appeared in \cite{GhysCircle}. See \cite[\S7.3]{GhysCircleSurvey} for an exposition.
\end{notes}

 \section{Additional ingredients of Ghys's proof}
 \label{GhysMissingSect}

\subsection{Amenability and an equivariant map}
\label{BdryMapSubsec}

A group is amenable if its action on every compact, convex set has a fixed point. More precisely:

\begin{defn}
A Lie group~$G$ is \defit{amenable} if, for every continuous action of~$G$ by linear operators on a locally convex topological vector space~$\TVS$, and every nonempty, compact, convex, $G$-invariant subset~$C$ of~$\TVS$,  the group~$G$ has a fixed point in~$C$.
\end{defn}

\begin{eg} \label{AmenEgs} \ 
\noprelistbreak
\begin{enumerate}
\item If $T$ is any continuous linear operator on~$\TVS$, and $v$~is any element of~$\TVS$, such that $\{T^n v\}$ is bounded, then every accumulation point of the sequence
	$$ v_n = \frac{1}{n} (Tv + T^2 v + \cdots + T^n v) $$
is a fixed point for~$T$. This implies that cyclic groups are amenable.
\item A generalization of this argument shows that all abelian groups are amenable; this statement is a version of the classical Kakutani-Markov Fixed-Point Theorem.
\item It is not difficult to see that if $N$ is a normal subgroup of~$G$, such that $N$ and $G/N$ are both amenable, then $G$ is amenable.
\item Combining the preceding two observations implies that solvable groups are amenable.
\item \label{AmenEgs-P}
In particular, the group~$P$ of \cref{SL3Notation} is amenable.
\end{enumerate}
\end{eg}

\begin{proof}[\normalfont\textbf{Proof of \cref{BdryMap}}]
Since $\Prob(X)$ is a compact, convex set, a version of the Banach-Alaoglu Theorem tells us that $L^\infty(G; \Prob(X) \bigr)$ is compact and convex in a natural weak$^*$ topology.
Let
	$$ L^\infty_\Gamma \bigl(G; \Prob(X) \bigr) 
		= \bigset{ \psi \in L^\infty(G; \Prob(X) \bigr) }{ 
		 \text{$\psi$ is $\Gamma$-equivariant (a.e.)} 
		} .$$
This is a closed, subset of $L^\infty(G; \Prob(X) \bigr)$, so it is compact. It is also convex and nonempty. To say $\psi$ is $\Gamma$-equivariant (a.e.) means, for each $\gamma \in \Gamma$, that $\psi(\gamma x) = \gamma \cdot \psi(x)$ for a.e.\ $x \in G$; so $G$ acts on $L^\infty_\Gamma \bigl(G; \Prob(X) \bigr)$ by translation on the right.
Hence, the subgroup~$P$ acts on $L^\infty_\Gamma \bigl(G; \Prob(X) \bigr)$.

Since $P$ is amenable (see \fullcref{AmenEgs}{P}), it must have a fixed point~$\psi_0$ in the compact, convex set $L^\infty_\Gamma \bigl(G; \Prob(X) \bigr)$. Then $\psi_0$ is invariant (a.e.) under translation on the right by elements of~$P$, so it factors through (a.e.)
to a well-defined map $\psi \colon G/P \to
\Prob(X)$. Because $\psi_0$ is $\Gamma$-equivariant, it is
immediate that $\psi$~is $\Gamma$-equivariant.
 \end{proof}

\subsection{Moore Ergodicity Theorem}
\label{MooreErgPfSect}

We will obtain the Moore Ergodicity Theorem \pref{MooreErgThm} as an easy consequence of the following result in representation theory:

\begin{thm}[(Decay of matrix coefficients)] \label{DecayMatCoeffSLn}
 If
 \noprelistbreak
 \begin{itemize}
 \item $G = \SL(n,\real)$,
 \item $\pi$ is a unitary representation of~$G$ on a Hilbert
space~$\hilbert$, such that no nonzero vector is fixed by $\pi(G)$; and
 \item $\{g_j\}$ is a sequence of elements of~$G$, such that $\lVert g_j
\rVert \to \infty$,
 \end{itemize}
 then $\langle \pi(g_j) \phi \mid \psi \rangle \to 0$, for every $\phi,\psi
\in \hilbert$.
 \end{thm}

\begin{proof}
 By passing to a subsequence, we may assume $\pi(g_j)$ converges weakly, to
some operator~$E$; that is,
 $$ \langle \pi(g_j) \phi \mid \psi \rangle
 \to \langle E \phi \mid \psi \rangle
 \mbox{ \ for every $\phi,\psi \in \hilbert$} .$$
We wish to show $\ker E = \hilbert$.

 Let 
 \begin{align*} \label{horodefn}
 U &= \{\, u \in G \mid g_j u g_j^{-1} \to e \,\} 
 \intertext{and} 
 U^- &= \{\, v \in G \mid g_j^{-1} v g_j \to e \,\} 
. \end{align*}
 For $u \in U$, we have
 \begin{align*}
 \langle E\pi(u) \phi \mid \psi \rangle
 &= \lim \langle \pi(g_j u) \phi \mid \psi \rangle
 \\&= \lim \langle \pi(g_j u g_j^{-1}) \pi(g_j) \phi \mid \psi \rangle
 \\&= \lim \langle \pi(g_j) \phi \mid \psi \rangle
 \\&= \langle E \phi \mid \psi \rangle 
 , \end{align*}
 so $E \pi(u) = E$. Therefore, letting $\hilbert^{U}$ be the space of $U$-invariant vectors in~$\hilbert$, we have
 	$$(\hilbert^{U})^\perp \subset \ker E .$$

We have
 $$ \langle E^* \phi \mid \psi \rangle
 =  \langle \phi \mid E \psi \rangle
 = \lim \langle \phi \mid \pi(g_j) \psi \rangle
 = \lim \langle \pi(g_j^{-1}) \phi \mid  \psi \rangle ,$$
 so the same argument, with $E^*$ in the place of~$E$ and $g_j^{-1}$ in the
place of~$g_j$, shows that 
	$$(\hilbert^{U^-})^\perp \subset \ker E^* .$$

Assume, for simplicity, that each $g_j$ is a positive-definite diagonal matrix:
	$$ g_j = \begin{bmatrix} a_j && \\ &b_j & \\ &&c_j \end{bmatrix} 
	\qquad \text{with $a_j,b_j,c_j > 0$}.$$
(It is not difficult to eliminate this hypothesis, by using the \emph{Cartan decomposition} $G = KAK$, but that is not necessary for Ghys's proof.) 
Then the subgroup generated by $\{\pi(g_j)\}$ is commutative.
 Because $\pi$ is unitary, this means that $\pi(g_j)$ commutes with both $\pi(g_k)$ and $\pi(g_k)^* = \pi(g_k^{-1})$ for every~$j$ and~$k$. 
Therefore, the limit~$E$ commutes with its adjoint (that is, $E$ is normal): we have $E^* E = E
E^*$. Hence
 \begin{align*}
  \lVert E\phi\rVert^2
 &= \langle E\phi \mid E \phi \rangle
 = \langle (E^* E)\phi \mid \phi \rangle
 \\&= \langle (E E^*)\phi \mid \phi \rangle
 = \langle E^*\phi \mid E^* \phi \rangle
 = \lVert E^* \phi\rVert^2 ,
 \end{align*}
 so $\ker E = \ker E^* $. 

Thus, 
 \begin{align*}
 \ker E
 &= \ker E + \ker E^*
 \\&\supset (\hilbert^{U})^\perp + (\hilbert^{U^-})^\perp
 \\&= (\hilbert^{U} \cap \hilbert^{U^-})^\perp
 \\&= (\hilbert^{\langle U, U^- \rangle})^\perp
 . \end{align*}

By passing to a subsequence, and then permuting the basis vectors of $\real^3$, we may assume 
	$$a_j \ge b_j \ge c_j .$$
Since $\|g_j\| \to \infty$, we have
	$$ \lim_{j \to \infty} \max \left\{  \frac{a_j}{b_j} ,  \frac{b_j}{c_j}  \right\} = \infty .$$
For definiteness, let us assume 
	$$ \text{$\displaystyle \limsup_{j \to \infty} \frac{a_j}{b_j} < \infty$ \quad and \quad  $\displaystyle \lim_{j \to \infty} \frac{b_j}{c_j} = \infty$},$$
 so
 	$$ U = \begin{bmatrix} 1&&\\ &1&\\ *&*&1 \end{bmatrix}
	\text{\qquad and \qquad}
	U^- = \begin{bmatrix} 1&&*\\ &1&*\\ &&1 \end{bmatrix} .$$
 (Other cases are similar.) Then it is easy to see that $\langle U, U^-
\rangle = G$, which means $\hilbert^{\langle U, U^- \rangle} =
\hilbert^G = 0$,
 so 
 	$$\ker E
	\supset (\hilbert^{\langle U, U^- \rangle})^\perp
 = 0^\perp
 = \hilbert ,$$
  as desired.
 \end{proof}
 
 \begin{proof}[\normalfont\textbf{Proof of \cref{MooreErgThm}}] \label{MooreErgThmPf}
 Suppose there is a $\Gamma$-invariant, measurable function on $G/H$ that is not constant (a.e.). Then:
 	\begin{align*}
	& \text{\hphantom{so }$\exists$ measurable function on $\Gamma \backslash G/H$ that is not constant (a.e.),}
	\\ & \text{so $\exists$ measurable function on $H \backslash G/\Gamma$ that is not constant (a.e.),}
	\\& \text{so $\exists$ $H$-invariant, measurable function~$f$ on $G/\Gamma$ that is not constant (a.e.).}
	\end{align*}
There is no harm in assuming that $f$ is bounded. Since $\Gamma$ is a lattice in~$G$, we know $G/\Gamma$ has finite measure, so this implies $f \in L^2(G/\Gamma)$. Letting $\pi$ be the natural unitary representation of~$G$ on $L^2(G/\Gamma)$, we know that $f$ is $\pi(H)$-invariant.

Let $L^2(G/\Gamma)_0$ be the orthogonal complement of the constant functions. Since $f$ is nonconstant, its projection~$\overline f$ in $L^2(G/\Gamma)_0$ is nonzero. By normalizing, we may assume $\| \overline f\|  = 1$. Since the orthogonal projection commutes with every unitary operator that preserves the space of constant functions, we know that $\overline f$, like~$f$, is $\pi(H)$-invariant. So
	$$ \text{$\langle \pi(h_j) \overline f \mid \overline f \rangle = \langle \overline f  \mid \overline f \rangle= 1$, \quad for all $h_j \in H$.} $$
On the other hand, since $H$~is closed and noncompact, we may choose a sequence $\{h_j\}$ of elements of~$H$, such that $\|h_j\| \to \infty$. Then, since no nonzero vector in  $L^2(G/\Gamma)_0$ is fixed by $\pi(G)$, \cref{DecayMatCoeffSLn} tells us that 
	$$ \langle \pi(h_j) \overline f \mid \overline f \rangle \to 0 .$$
This is a contradiction.
\end{proof}

\begin{rem}
The assumption that $H$ is not compact is necessary in the Moore Ergodicity Theorem \pref{MooreErgThm}: it is not difficult to see that if $H$ is a \emph{compact} subgroup of $G = \SL(3,\real)$, and $\Gamma$ is  any lattice in~$G$, then $\Gamma$ is \emph{not} ergodic on $G/H$.
\end{rem}

\begin{cor} \label{EquiPsiAtoms}
In the situation of \cref{BdryMap}, one may assume that either:
 \noprelistbreak
 \begin{itemize}
 \item for a.e.\ $x \in G/P$, the measure $\psi_0(x)$ has no atoms, 
 or
 \item for a.e.\ $x \in G/P$, the measure $\psi_0(x)$ consists
entirely of atoms.
 \end{itemize}
\end{cor}

\begin{proof}
For each $x \in G/P$, write $\psi(x) = \psi_{\text{noatom}}(x)
+ \psi_{\text{atom}}(x)$, where $\psi_{\text{noatom}}(x)$ has no
atoms, and $\psi_{\text{atom}}(x)$~consists entirely of atoms. 
Since $\psi_{\text{noatom}}$ and~$\psi_{\text{atom}}(x)$ are uniquely determined by~$\psi$, it is not difficult to see that they, like~$\psi$, are measurable and $\Gamma$-equivariant. One or the other must be nonzero on a set of positive measure, and then the ergodicity of $\Gamma$ on $G/P$ implies that this function is nonzero almost everywhere, so it can be normalized to define a ($\Gamma$-equivariant, measurable) map into $\Prob(\circle)$.
\end{proof}

\begin{cor} \label{ErgOn(G/P)2}
In the situation of \cref{SL3Notation},\/ $\Gamma$ is ergodic on $(G/P)^2$.
\end{cor}

\begin{proof}
Two flags $F_1 = (\ell_1,\pi_1)$ and $F_2 = (\ell_2,\pi_2)$ are in \defit{general position} if $\ell_1 \notin \pi_2$ and $\ell_2 \notin \pi_1$. It is not difficult to see that $G$ is transitive on the set~$\flags^2_0$ of pairs of flags in general position, and that the complement of~$\flags^2_0$ has measure~$0$ in $\flags^2$. Therefore, $(G/P)^2$ may be identified (a.e.) with $G/H$, where $H$ is the stabilizer of some pair of flags in general position; we may take
	$$ H 
	= \Stab_G \Bigl( \bigl( (*,0,0), (*,*,0) \bigr), \bigl( (0,0,*), (0,*,*) \bigr) \Bigr)
	= \begin{bmatrix} *&&\\ &*&\\ &&* \end{bmatrix} .$$
Since $H$ is not compact, the Moore Ergodicity Theorem \pref{MooreErgThm} tells us that $\Gamma$ is ergodic on $G/H \approx (G/P)^2$.
\end{proof}


\begin{notes}
\Cref{BdryMap} is due to Furstenberg \cite{Furstenberg-Bdry}.
It is a basic result in the theory of lattices, so proofs can be found in numerous references, including \cite[Prop.~7.11]{GhysCircleSurvey} and \cite[Prop.~4.3.9, p.~81]{ZimmerBook}.

The monograph \cite{PierBook} is a standard reference on amenability. See \cite[\S4.1]{ZimmerBook} for a brief treatment.

\Cref{MooreErgThm} is due to C.\,C.\,Moore \cite{Moore-ergodicity}. 
The stronger \cref{DecayMatCoeffSLn} is due to R.\,Howe and C.\,C.\,Moore \cite[Thm.~5.1]{HoweMoore} and (independently) R.\,J.\,Zimmer
\cite[Thm.~5.2]{Zimmer-orbitspace}. The elementary proof we give here was found by R.\,Ellis and M.\,Nerurkar \cite{EllisNerurkar}. 
\end{notes}

\section{Bounded cohomology and the Burger-Monod proof} \label{BurgerMonodPf}

\subsection{Bounded cohomology and actions on the circle}

Suppose a discrete group~$\Gamma$ acts by orientation-preserving homeomorphisms on $\circle = \real/\integer$. Since $\real$ is the universal cover of~$\circle$, each element~$\gamma$ of~$\Gamma$ can be lifted to a homeomorphism~$\widetilde\gamma$ of~$\real$. The lift~$\widetilde\gamma$ is not unique, but it is well-defined if we require that $\widetilde\gamma(0) \in [0,1)$. 

\begin{defn} \label{EulerCocycDefn}
For $\gamma_1,\gamma_2 \in \Gamma$, the homeomorphisms $\widetilde{\gamma_1 \gamma_2}$ and $\widetilde{\gamma_1} \widetilde{\gamma_2}$ are lifts of the same element $\gamma_1 \gamma_2$ of~$\Gamma$, so there exists $z = z(\gamma_1, \gamma_2) \in \integer$, such that 
	$$\widetilde{\gamma_1} \widetilde{\gamma_2} = \widetilde{\gamma_1 \gamma_2}  + z .$$
The map $z \colon \Gamma \times \Gamma \to \integer$ is called the \defit{Euler cocycle} of the action of~$\Gamma$ on~$\circle$.
\end{defn}

It is easy to see that 
	\noprelistbreak
	\begin{itemize}
	\item $z$ is a bounded function (in fact, $z(\Gamma \times \Gamma) \subset \{0,1\}$),
	and
	\item $z(\gamma_1,\gamma_2) + z(\gamma_1\gamma_2,\gamma_3)
	=  z(\gamma_1, \gamma_2 \gamma_3) + z(\gamma_2,\gamma_3)$, so $z$~is an Eilenberg-MacLane $2$-cocyle.
	\end{itemize}
Therefore, the Euler cocycle determines a bounded cohomology class:

\begin{defn} \ 
	\noprelistbreak
	\begin{enumerate}
	\item A \defit{bounded $k$-cochain} is a bounded function $c \colon \Gamma^k \to \integer$.
	\item The bounded cochains form a chain complex with respect to the differential 
		$$\delta \colon \Cbdd^k(\Gamma;\integer) \to  \Cbdd^{k+1}(\Gamma;\integer) $$  defined by
		\begin{align*}
		 \delta c (\gamma_0,\gamma_1,\ldots,\gamma_k) 
		&= 
		c (\gamma_1,\ldots,\gamma_k) 
		\\& \qquad + \sum_{i=1}^{k} (-1)^{i} c(\gamma_0,\ldots,\gamma_{i-1} \gamma_{i},\ldots,\gamma_k)
		\\& \qquad \qquad + (-1)^{k+1} c (\gamma_0,\gamma_1,\ldots,\gamma_{k-1}) 
		. \end{align*}
The cohomology of this complex is the \defit{bounded cohomology} of~$\Gamma$.
	\item The \defit{Euler class} of the action of~$\Gamma$ is the cohomology class $[z] \in \Hbdd^2(\Gamma,\integer)$ determined by the Euler cocycle.
	\end{enumerate}
\end{defn}

\begin{rem}
The Euler cocycle~$z$ depends on the choice of the covering map from~$\real$ to~$\circle$, but it is not difficult to see that the Euler class $[z]$ is well-defined. Indeed, it is an invariant of the (orientation-preserving) homeomorphism class of the action.
\end{rem}

The connection with \cref{GhysBurgerMonodFinOrb} is provided by the following fundamental observation:

\begin{prop}[(Ghys)] \label{BddCoho<>FP}
The action of\/~$\Gamma$ on~$\circle$ has a fixed point if and only if its Euler class is $0$ in $\Hbdd^2(\Gamma;\integer)$.
\end{prop}

\begin{proof}
($\Rightarrow$)
We may assume the fixed point is the image of~$0$ under the covering map $\real \to \circle$. Then $\widetilde \gamma(0) = 0$, for every $\gamma \in \Gamma$, so it is clear that $z(\gamma_1,\gamma_2) = 0$ for all $\gamma_1$ and~$\gamma_2$.

($\Leftarrow$) Assume $z = \delta \varphi$, where $\varphi \colon \Gamma \to \integer$ is bounded. If we set $\widehat\gamma = \widetilde\gamma - \varphi(\gamma)$, then the map $\gamma \to \widehat\gamma$ is a homomorphism; $\widehat\Gamma$ is a lift of~$\Gamma$ to a group of homeomorphisms of~$\real$. 

Since $\widetilde\gamma(0) \in [0,1)$, and $\varphi$ is bounded, it is clear that the $\widehat\Gamma$-orbit of~$0$ is bounded. 
Since the orbit is obviously a $\widehat\Gamma$-invariant set, its supremum is also $\widehat\Gamma$-invariant; in other words, the supremum is a fixed point for~$\widehat\Gamma$ in~$\real$. The image of this fixed point under the covering map is a fixed point for~$\Gamma$ in~$\circle$.
\end{proof}

The definition of $\Hbdd^k(\Gamma;\integer)$ can be generalized to allow any coefficient module in the place of~$\integer$. For real coefficients, we have the following important fact:

\begin{cor} \label{Hbdd=0->FinOrb}
If $\Hbdd^2(\Gamma;\real) = 0$, and the abelianization of\/~$\Gamma$ is finite, then every action of\/~$\Gamma$ on\/~$\circle$ has a finite orbit.
\end{cor}

\begin{proof}
The short exact sequence
	$$ 0 \to \integer \to \real \to \real/\integer \to 0$$
of coefficient groups leads to a long exact sequence of bounded cohomology groups. A part of this sequence is
	$$ \Hbdd^1(\Gamma; \real/\integer) \to \Hbdd^2(\Gamma;\integer) \to \Hbdd^2(\Gamma;\real) .$$
The right end of this sequence is~$0$, by assumption. 
If we assume, for simplicity, that the abelianization of~$\Gamma$ is trivial (rather than merely finite), then the left end is also~$0$. Hence, the middle term must be~$0$. Then \cref{BddCoho<>FP} implies that every (orientation-preserving) action of~$\Gamma$ on~$\circle$ has a fixed point.

Without the simplifying assumption, one can obtain the weaker conclusion that the commutator subgroup $[\Gamma,\Gamma]$ has a fixed point. Since, by hypothesis, $\Gamma/[\Gamma,\Gamma]$ is finite, this implies that the action of~$\Gamma$ has a finite orbit.
\end{proof}


We now need two observations:
	\noprelistbreak
	\begin{enumerate}
	\item Forgetting that a bounded $k$-cocycle is bounded yields a natural map from bounded cohomology to ordinary cohomology:
	\begin{equation} \label{BddCohoComparison}
	\tau_\Gamma \colon \Hbdd^k(\Gamma;\real) \to H^k(\Gamma;\real) 
	. \end{equation}
	\item The cohomology of lattices in $\SL(n,\real)$ has been studied extensively; in particular, it is known that
	\begin{equation} \label{H2(LattSLn;R)=0}
	 \text{$H^2(\Gamma;\real) = 0$ if $\Gamma$ is any lattice in $\SL(n,\real)$, with $n \ge 6$.}
	 \end{equation}
	\end{enumerate}
Therefore, under the assumption that $n \ge 6$,
the conclusion of \cref{GhysBurgerMonodFinOrb} can be obtained by combining \cref{Hbdd=0->FinOrb} with the following result.

\begin{thm}[(Burger-Monod)] \label{BurgerMonodInjectSLn}
If\/ $\Gamma$ is any lattice in\/ $\SL(n,\real)$, with $n \ge 3$, then the comparison map\/ \pref{BddCohoComparison} is injective for $k = 2$.
\end{thm}

\begin{cor}[(Burger-Monod)] \label{BurgerMonodVanishSLn}
If\/ $\Gamma$ is any lattice in\/ $\SL(n,\real)$, with $n \ge 3$, then $\Hbdd^2(\Gamma;\real) = 0$.
\end{cor}

\begin{rem} \label{BMWithoutVanishing}
See \fullcref{GhysBurgerMonodGeneralRem}{criterium} for a brief mention of how to obtain \cref{GhysBurgerMonodFinOrb} from \cref{BurgerMonodInjectSLn}, without needing to know that $H^2(\Gamma;\real)$ vanishes.
\end{rem}

\subsection{Outline of the Burger-Monod proof of injectivity}

M.\,Burger and N.\,Monod developed an extensive and powerful theory for the study of bounded cohomology, but we will discuss only the parts that are used in the proof of \cref{BurgerMonodInjectSLn}, and even these will only be sketched.

\begin{assump} In the remainder of this section:
	\noprelistbreak
	\begin{itemize}
	\item $G = \SL(n,\real)$, with $n \ge 3$,
	and
	\item $\Gamma$ is a lattice in~$G$.
	\end{itemize}
To avoid a serious technical complication, we will assume that $G/\Gamma$ is compact.
\end{assump}

\begin{proof}[Outline of the proof of \cref{BurgerMonodInjectSLn}]
We employ relations between the cohomology of~$\Gamma$ and the cohomology of~$G$. (The bounded cohomology of~$G$ will be introduced in \cref{BddCohoGDefn} below. When working with $G$, we always use \emph{continuous} cochains.)
	\noprelistbreak
	\begin{itemize}
	\item We will see that there is a natural map $i_{\text{\normalfont\upshape bdd}} \colon \Hbdd^k(\Gamma;\real) \to \Hbdd^k \bigl( G; L^2(G/\Gamma) \bigr)$.
	\item It is a classical fact that if $G/\Gamma$ is compact, then there is a natural map $i \colon H^k(\Gamma;\real) \to H^k \bigl( G; L^2(G/\Gamma) \bigr)$.
	\item We have comparison maps $\tau_\Gamma$ and $\tau_G$ from bounded cohomology to ordinary cohomology.
	\end{itemize}
Letting $k = 2$ yields the following commutative diagram:
	$$\begin{CD}
	\Hbdd^2(\Gamma;\real) @>i_{\text{\normalfont\upshape bdd}}>> \Hbdd^2 \bigl( G; L^2(G/\Gamma) \bigr) \\
	@V\tau_\Gamma VV     @VV \tau_G V     \\
	H^2(\Gamma;\real) @>i>> H^2 \bigl( G; L^2(G/\Gamma) \bigr)
	\end{CD}$$
We will show that $i_{\text{\normalfont\upshape bdd}}$ and $\tau_G$ are both injective \csee{L2IndInj,CompareGInj}. The commutativity of the diagram then implies that $\tau_\Gamma$ is also injective.
\end{proof}

\subsubsection*{Injectivity of $i_{\text{\normalfont\upshape bdd}}$}
Cohomology, whether bounded or not, can be described either in terms of inhomogeneous cocycles, or in terms of homogeneous cocycles. The Euler cocycle arose in \cref{EulerCocycDefn} as an inhomogeneous cocycle, but the injectivity of~$i_{\text{\normalfont\upshape bdd}}$ is easier to explain in homogeneous terms. The following definition is written with real coefficients, because we no longer have any need for $\integer$-coefficients in our discussion.

\begin{defn} \label{HomogBddCohoDefn} \ 
	\noprelistbreak
	\begin{enumerate}
	\item A \defit{homogeneous bounded $k$-cochain} on~$\Gamma$ is a bounded function $\dot c \colon \Gamma^{k+1} \to \real$, such that 
		$$ \dot c ( \gamma \gamma_0, \gamma \gamma_1, \ldots, \gamma \gamma_{k} ) 
		= \dot c ( \gamma_0, \gamma_1, \ldots, \gamma_{k} ) ,$$
	for all $\gamma, \gamma_0, \gamma_1, \ldots, \gamma_{k} \in \Gamma$.
	\item The homogeneous bounded cochains form a chain complex with respect to the differential 
		$$\delta \colon \hCbdd^k(\Gamma; \real) \to  \hCbdd^{k+1}(\Gamma; \real) $$  defined by
		\begin{align*}
		 \delta \dot c (\gamma_0,\gamma_1,\ldots,\gamma_{k+1}) 
		= 
		\sum_{i=0}^{k+1} (-1)^{i} \dot c(\gamma_0,\ldots, \widehat{\gamma_{i}},\ldots,\gamma_{k+1})
		, \end{align*}
where $\widehat{\gamma_{i}}$ denotes that $\gamma_i$ is omitted.
	\end{enumerate}
For any bounded $k$-cochain~$c$, there is a corresponding homogeneous bounded $k$-cochain $\dot c$, defined by
	$$ \dot c(\gamma_0, \gamma_1, \ldots, \gamma_{k} ) 
	= c( \gamma_0^{-1} \gamma_1, \gamma_1^{-1} \gamma_2, \ldots, \gamma_{k-1}^{-1} \gamma_k) .$$
Thus, the cohomology of the complex $\{\hCbdd^k(\Gamma; \real)\}$ is the {bounded cohomology} of~$\Gamma$.
\end{defn}

\begin{notation} \ 
\noprelistbreak
\begin{enumerate}
\item  $P$ is the group of upper-triangular matrices in $G = \SL(n,\real)$ \ccf{SL3Notation}.
\item $ZL^\infty_{\text{\normalfont\upshape alt}} \bigl( (G/P)^3; \real \bigr)^\Gamma$ is the vector space of all $f \in L^\infty \bigl( (G/P)^3; \real \bigr)$, such that
	\noprelistbreak
	\begin{enumerate}
	\item $f$ is \emph{alternating}; i.e., 
		$$f( x_{\sigma(1)},  x_{\sigma(2)},  x_{\sigma(3)}) = \mathrm{sgn}(\sigma) \, f(x_1,x_2,x_3) $$
	for every permutation~$\sigma$ of $\{1,2,3\}$,
	\item $f$ is \emph{$\Gamma$-invariant}; i.e., for every $\gamma \in \Gamma$, we have 
		$$f(\gamma x_1, \gamma x_2, \gamma x_3) = f(x_1,x_2,x_3)$$
		for a.e.\ $x_1,x_2,x_3$,
		and
	\item $f$ is a \emph{cocycle}; i.e., for a.e.\ $x_0,x_1,x_2,x_3$, we have
		$$ \sum_{i=0}^3 (-1)^i f(x_0,\ldots, \widehat{x_i}, \ldots,x_3) = 0 .$$
	\end{enumerate}
\end{enumerate}
\end{notation}

\begin{rem}
If $\Gamma$ were ergodic on $(G/P)^3$, then the following theorem would immediately imply that $\Hbdd^2(\Gamma;\real) = 0$. This is because any $\Gamma$-invariant function on $(G/P)^3$ would have to be constant, so could not be alternating (unless it were identically~$0$). 
\end{rem}

\begin{thm}[(Burger-Monod)] \label{H2bdd=ZLinfty}
$\Hbdd^2(\Gamma;\real) \iso ZL^\infty_{\text{\normalfont\upshape alt}} \bigl( (G/P)^3; \real \bigr)^\Gamma$.
\end{thm}

\begin{proof}
For each homogeneous bounded $k$-cocycle $\dot c \colon \Gamma^{k+1} \to \real$, we will show how to construct a corresponding $\check c \in ZL^\infty_{\text{\normalfont\upshape alt}}\bigl( (G/P)^{k+1}; \real \bigr)^\Gamma$. The map $\dot c \mapsto \check c$ intertwines the differentials, so it induces a map from $\Hbdd^k(\Gamma;\real)$ to the cohomology of the chain complex 
	$$ \left\{\, ZL^\infty_{\text{\normalfont\upshape alt}}\bigl( (G/P)^{k+1}; \real \bigr)^\Gamma \,\right\} ,$$
and, although we will not prove it, this map is an isomorphism on cohomology.

Since $\Gamma$ is ergodic on $(G/P)^2$ \csee{ErgOn(G/P)2}, every $\Gamma$-invariant function on~$(G/P)^2$ is constant, so $0$ is the only such function that is alternating. Hence, there are no coboundaries in $ZL^\infty_{\text{\normalfont\upshape alt}}\bigl( (G/P)^{3}; \real \bigr)^\Gamma$. This establishes the conclusion of the theorem.

\medbreak

To complete the proof, we now describe the construction of~$\check c$. For simplicity, let us assume $k = 2$, so $\dot c \colon \Gamma^3 \to \real$ is a homogeneous bounded $2$-cocycle. By making use of the cocycle identity
	$$ \dot c(x_1,x_2,x_3) = \dot c(x_0,x_2,x_3) - \dot c(x_0,x_1,x_3) + \dot c(x_0,x_1,x_2) ,$$
one can show that $\dot c$ is alternating.
Therefore, $\dot c$ can be extended to 
	$$\bar c \in ZL^\infty_{\text{\normalfont\upshape alt}}( G^3; \real)^\Gamma ,$$
by choosing a fundamental domain~$\fund$ for~$\Gamma$ in~$G$ and making $\bar c$ constant on $\gamma_1 \fund \times \gamma_2 \fund \times \gamma_3 \fund $, for all $\gamma_1,\gamma_2,\gamma_3 \in \Gamma$.

Now, because $P$ is amenable \fullcsee{AmenEgs}{P}, there is a left-invariant mean~$\mu$ on $L^\infty(P)$. Using $\mu$ to average on each left coset of~$P$ yields a map 
	$$\overline \mu \colon L^\infty(G) \to L^\infty(G/P) .$$
Then a map
	$$\overline \mu^3 \colon L^\infty(G^3) \to L^\infty\bigl( (G/P)^3 \bigr) $$
can be constructed by, roughly speaking, setting $\overline \mu^3 = \overline \mu \otimes \overline \mu \otimes \overline \mu$. Let $\check c = \overline \mu^3 (\bar c) \in L^\infty \bigl( (G/P)^3 ; \real \bigr)$.
\end{proof}

\begin{defn} \label{BddCohoGDefn} \ 
\noprelistbreak
\begin{enumerate}
\item The notion of a \emph{homogeneous bounded $k$-cochain} on~$G$ is defined by replacing $\Gamma$ with~$G$ in \cref{HomogBddCohoDefn}, and requiring $\dot c$ to be \emph{continuous}.
\item The cohomology of the complex $\hCbdd^k(G;\real)$ is $\Hbdd(G;\real)$, the\/ \emph{{\upshape(}continuous\/{\upshape)} bounded cohomology} of~$G$.
\end{enumerate}
\end{defn}

\begin{cor}[(Burger-Monod)] \label{L2IndInj}
There is a natural injection 
	$$ \Hbdd^2(\Gamma; \real) \hookrightarrow \Hbdd^2 \bigl( G; L^2(G/\Gamma) \bigr) .$$
\end{cor}

\begin{proof}
For any $\dot c \in ZL^\infty_{\text{\normalfont\upshape alt}}\bigl( (G/P)^{3}; \real \bigr)^\Gamma$ and $x \in (G/P)^3$, we can define 
	$$ \text{$\dot c_x \in L^\infty(G/\Gamma) \subset L^2(G/\Gamma)$ 
	\quad by \quad 
	$\dot c_x(g\Gamma) = \dot c(gx)$.}$$
The map $x \mapsto \dot c_x$ is $G$-equivariant, so it is an element of $ZL^\infty_{\text{\normalfont\upshape alt}}\bigl( (G/P)^{3}; L^2(G/\Gamma) \bigr)^G$. Therefore, we have an injection 
	$$ ZL^\infty_{\text{\normalfont\upshape alt}}\bigl( (G/P)^{3}; \real \bigr)^\Gamma
	\hookrightarrow 
	ZL^\infty_{\text{\normalfont\upshape alt}}\bigl( (G/P)^{3}; L^2(G/\Gamma) \bigr)^G
	. $$
\Cref{H2bdd=ZLinfty} identifies the domain of this injection with $\Hbdd(\Gamma;\real)$, and the same argument identifies the target with $\Hbdd\bigl( G; L^2(G/\Gamma) \bigr)$.
\end{proof}

\begin{thm}[(Burger-Monod)] \label{CompareGInj}
The comparison map
	$$ \tau_G \colon \Hbdd^2 \bigl( G; L^2(G/\Gamma) \bigr) \to H^2 \bigl( G; L^2(G/\Gamma) \bigr)$$
is injective.
\end{thm}

\subsubsection*{Injectivity of $\tau_G$}

The Hilbert space $L^2(G/\Gamma)$ decomposes as the direct sum of the constant functions~$\complex$ and the space $L^2_0(G/\Gamma)$ of functions with integral~$0$. The theorem is proved for the two summands individually \csee{CompareGInjC,CompareGInjL0}. In both cases, we will argue with inhomogeneous cochains.

\begin{prop} \label{CompareGInjC}
The comparison map
	$$ \Hbdd^2 \bigl( G; \complex \bigr) \to H^2 \bigl( G; \complex \bigr)$$
is injective.
\end{prop}

\begin{proof}
Let $c$ be an inhomogeneous cocycle that represents a class in the kernel of the comparison map. This implies that 
	\noprelistbreak
	\begin{itemize}
	\item $c \colon G \times G \to \complex$ is a bounded, continuous function,
	and
	\item there a a continuous function $\varphi \colon G \to \complex$, such that $\delta \varphi = c$.
	\end{itemize}
It suffices to show $\varphi$ is bounded, for then $c$~is the coboundary of a bounded cochain (namely,~$\varphi$), so $[c] = 0$ in bounded cohomology.

Note that, for all $g,h \in G$, we have
	\begin{equation} \label{quasimorphism}
	 |\varphi( gh ) - \varphi(g) - \varphi(h)| = |\delta \varphi(g,h)| =  |c(g,h)| \le \|c\|_\infty 
	 . \end{equation}

Now assume, for concreteness, that $G = \SL(3,\real)$, and let
	\begin{align*}
	 U_{1,2} = \begin{bmatrix} 1&*& \\ &1& \\ &&1 \end{bmatrix},
 	\quad
	U_{1,3} = \begin{bmatrix} 1&&* \\ &1& \\ &&1 \end{bmatrix},
  	\quad
	U_{2,1} = \begin{bmatrix} 1&& \\ *&1& \\ &&1 \end{bmatrix},
	\\
 	U_{2,3} = \begin{bmatrix} 1&& \\ &1&* \\ &&1 \end{bmatrix},
  	\quad
	U_{3,1} = \begin{bmatrix} 1&& \\ &1& \\ *&&1 \end{bmatrix},
  	\quad
	U_{3,2} = \begin{bmatrix} 1&& \\ &1& \\ &*&1 \end{bmatrix}
	. \end{align*}
We will show that $\varphi$ is bounded on~$U_{1,2}$, and a similar argument shows that $\varphi$ is bounded on each of the other elementary unipotent subgroups $U_{i,j}$. 
Then the desired conclusion that $\varphi$ is bounded on all of~$G$ is obtained by combining these bounds with \pref{quasimorphism} and the elementary observation that, for some $N \in \natural$, there exist $i_1,\ldots,i_N$ and $j_1,\ldots,j_N$, such that
	$$ G = U_{i_1,j_1} U_{i_2,j_2} U_{i_3,j_3} \cdots U_{i_N,j_N}  . $$

To complete the proof, we now show that $\varphi$ is bounded on~$U_{1,2}$. To see this, note that, for
	$$ a = \begin{bmatrix} 2 && \\ &1/2&\\ && 1 \end{bmatrix} , $$
	we have
	$$ \text{$\displaystyle \lim_{k \to \infty} a^{-k} u a^k = e$ \quad for all $u \in U_{1,2}$.} $$
Therefore, for $u \in U_{1,2}$ and $k \in \natural$, repeated application of \pref{quasimorphism} yields
	\begin{align*}
	|h(u)|
	& \le | h(a^{k}) +  h(a^{-k} u a^k) + h(a^{-k}) | + 2 \| c\|_\infty
	\\& \le | h(a^{-k} u a^k)| + | h(a^{k})  + h(a^{-k}) | + 2 \| c\|_\infty
	\\& \le | h(a^{-k} u a^k)| + | h(e) | + 3 \| c\|_\infty
	\\& \to |h(e)| + | h(e) | + 3 \| c\|_\infty && \text{as $k \to \infty$}
	.  \end{align*}
\end{proof}

\begin{thm}[(Burger-Monod)] \label{CompareGInjL0}
The comparison map
	$$  \Hbdd^2 \bigl( G; L^2_0(G/\Gamma) \bigr) \to H^2 \bigl( G; L^2_0(G/\Gamma) \bigr)$$
is injective.
\end{thm}

\begin{proof}
Let $c$ be an inhomogeneous cocycle that represents a class in the kernel of the comparison map, so $c = \delta \varphi$, for some continuous $\varphi \colon G \to L^2_0(G/\Gamma)$. It suffices to show that $\varphi$ is bounded.

Note that, letting $\pi$ be the representation of~$G$ on $L^2_0(G/\Gamma)$, we have, for all $g,h \in G$,
	\begin{equation} \label{CobdryBdd}
	 \| \, \varphi( gh ) - \varphi(g) - \pi(g) \, \varphi(h)| = |\delta \varphi(g,h)| =  |c(g,h) \,\| \le \|c\|_\infty 
	 . \end{equation}

Let us assume $G = \SL(5,\real)$. (The same argument works for all $n \ge 5$, but some modifications are needed when $n$~is small.) 
Much as in the proof of \cref{CompareGInjC}, it suffices to show that $\varphi$ is bounded on
	$$ U_{1,2} = \begin{bmatrix} 1&*& & & \\ &1& & & \\ &&1 & & \\ &&&1& \\ &&&&1 \end{bmatrix} .$$

	Let 
		$$ H =  \begin{bmatrix} 1&& & & \\ &1& & & \\ &&* &* &* \\ &&*&*&* \\ &&*&*&* \end{bmatrix} \iso \SL(3,\real) .$$
For all $u \in U_{1,2}$ and $h \in H$, we have
	\begin{align*}
	&\| \,\varphi(u) - \varphi(h) - \pi(h) \varphi(u) - \pi(h) \pi(u) \, \varphi(h^{-1}) \,\|
	\\& \qquad = 
	\bigl\| \,\bigl( \varphi(huh^{-1}) - \varphi(h)  - \pi(h) \, \varphi(uh^{-1})  \bigr)
	&& \begin{pmatrix}
	\text{$H$ commutes with~$U_{1,2}$,} \\ \text{so $u = h u h^{-1}$}
	\end{pmatrix}
	\\& \qquad \qquad {}
	+ \pi(h)  \, \bigl(\varphi(uh^{-1}) - \varphi(u) - \pi(u) \, \varphi(h^{-1}) \bigr) \,\bigr\|
	\\& \qquad \le 2 \|c\|_\infty
	.
	&& \text{(by \pref{CobdryBdd})}
	\end{align*}
So
	$$ \bigl\| \bigl( \Id - \pi(h) \bigr) \varphi(u) \bigr\| \le \| \varphi(h) \| + \| \varphi(h^{-1} \| + 2 \|c\|_\infty . $$
Hence, for any function~$f$ on~$H$ whose integral is~$1$, we have
	\begin{equation} \label{(Id-pi(h))phi(u)}
	 \bigl\| \bigl( \Id - \pi(f) \bigr) \, \varphi(u) \bigr\| \le 
	\sup_{h \in \supp f} \Bigl( \| \varphi(h) \| + \| \varphi(h^{-1} \| + 2 \|c\|_\infty \Bigr) 
	. \end{equation}

To complete the proof, we need to make a good choice of the function~$f$.
There are no $G$-invariant vectors in $L_0(G/\Gamma)$, so the Moore Ergodicity Theorem \pref{MooreErgThm} implies there are no $H$-invariant vectors (see the proof on page~\pageref{MooreErgThmPf}). Since $H = \SL(3,\real)$ has Kazhdan's property~$(T)$, we conclude that $H$ has no almost-invariant vectors. Hence, it is not difficult to see that there is a continuous function~$f$ on~$H$, such that
	\noprelistbreak
	\begin{itemize}
	\item $f$ has compact support, 
	\item $\int_H f \, d \mu_H = 1$ (where $\mu_H$ is the Haar measure on~$H$),
	and
	\item $\| \pi(f) \| < 1$.
	\end{itemize}
Then $\Id - \pi(f)$ is invertible, so \pref{(Id-pi(h))phi(u)} implies $\| \, \varphi(u) \,\|$ is bounded (independent of~$u$), as desired.
\end{proof}

\begin{notes}
See \cite{Monod-InviteBddCoho} for a recent introduction to bounded cohomology.
Although our presentation of bounded cohomology takes a very naive approach, the work of Ghys and Burger-Monod is much more functorial.

\Cref{BddCoho<>FP} is due to \'E.\,Ghys \cite{GhysBddCoho}. In fact, he proved that the Euler class determines the action up to semi-conjugacy.

The Burger-Monod theory of bounded cohomology (including \cref{BurgerMonodInjectSLn,GhysBurgerMonodFinOrb}) was developed in \cite{BurgerMonod-BddCoho1,BurgerMonod-BddCoho2}. An exposition appears in \cite{Monod-ContBddCoho}.
The improvement mentioned in \cref{BMWithoutVanishing} appears in the paper of M.\,Burger \cite{Burger-criterium} in this volume.

The cohomology vanishing result \pref{H2(LattSLn;R)=0} is a special case of \cite[Thm.~4.4(ii)]{Borel-StableCoho2}.

\end{notes}

\section[Bounded orbits and bounded generation]
{Non-orderability from bounded orbits and bounded generation}
\label{BddOrbBddGenSect}

Recall that \cref{Witte-AOCThm} can be stated in any one of the following three equivalent forms (cf.\ \cref{AOC<>LO,LO<>AOR}):
	\noprelistbreak
	\begin{enumerate} \renewcommand{\theenumi}{\Alph{enumi}}
	\item If $n \ge 3$, then no finite-index subgroup of $\SL(n,\integer)$ has a faithful action on~$\circle$.
	\item If $n \ge 3$, then no finite-index subgroup of $\SL(n,\integer)$ has an orientation-preserving, faithful action on~$\real$.
	\item If $n \ge 3$, then no finite-index subgroup of $\SL(n,\integer)$ is left orderable.
	\end{enumerate}
A short, elementary proof of this theorem was given in \cref{SLnZNotLOSect}, but we will now describe a different approach that has the potential to work in a more general situation. It has two main ingredients: 
	\noprelistbreak
	\begin{itemize}
	\item bounded orbits of unipotent elementary matrices (\cref{ElemOrbsBdd}),
	and
	\item bounded generation by unipotent elementary matrices (\cref{CKPThm}).
	\end{itemize}
To keep things simple, let us assume $n = 3$.

\begin{defn}
A matrix~$u$ in $\SL(3,\real)$ is a \defit{unipotent elementary matrix} if 
	\noprelistbreak
	\begin{itemize}
	\item $u_{1,1} = u_{2,2} = u_{3.3} = 1$ (i.e., $u$ has all $1$'s on the main diagonal), 
	and
	\item $u$ has only one nonzero entry off the main diagonal.
	\end{itemize}
(In other words, a unipotent elementary matrix is one of the matrices $a_1,\ldots,a_6$ of \pref{ElemUnipMatsInSL3}, for some $k \in \real$.)
\end{defn}

\begin{prop} \label{ElemOrbsBdd}
Suppose a finite-index subgroup\/~$\Gamma$ of\/ $\SL(3,\integer)$ has an orientation-preserving, faithful action on~$\real$.

If $u$ is any unipotent elementary matrix in~$\Gamma$, then the $u$-orbit of each point in~$\real$ is a bounded set.
\end{prop}

\begin{proof}
It suffices to show that the $u$-orbit of~$0$ is bounded above.

Define a left-invariant total order~$\prec$ on~$\Gamma$, as in \cref{LO<>AOR}($\AORimpliesLO$), so
	$$ a(0) < b(0) \quad \Rightarrow \quad a \prec b .$$
By permuting the standard basis vectors, we may assume
	$$ u = \begin{bmatrix} 1&&*\\ &1\\ &&1 \end{bmatrix} .$$
Then \cref{LOHeis} implies there is some $a \in \Gamma$, such that $u \ll a$ (see \pref{SLnZnotLOPf-<<} with $i = 2$); i.e., $u^k \prec a$, for all $k \in \integer$. From the definition of~$\prec$, this means that
	$$ \text{$u^k(0) < a(0)$, \quad for all~$k \in \integer$} .$$
So the $u$-orbit of~$0$ is bounded above (by $a(0)$).
\end{proof}

It is a standard fact of undergraduate linear algebra that any invertible matrix is a product of elementary matrices. (This is a restatement of the fact that every invertible matrix can be reduced to the identity by elementary row operations.) Furthermore, it is not difficult to see that  if the invertible matrix has determinant~$1$, then:
	\noprelistbreak
	\begin{itemize}
	\item the elementary matrices can be assumed to be unipotent,
	and
	\item there is a bound on the  \emph{number} of elementary matrices that is needed.
	\end{itemize}

\begin{exer}
Show that every matrix in $\SL(3,\real)$ is a product of $< 10$ unipotent elementary matrices.
\end{exer}

It is a much deeper fact that the boundedness remains true when the field~$\real$ is replaced with the ring~$\integer$:

\begin{thm}[(Carter-Keller)] \label{SL3ZBddGen}
Every element of\/ $\SL(3,\integer)$ is a product of $< 50$ unipotent elementary matrices in~$\SL(3,\integer)$.
\end{thm}

There is also a bound for any finite-index subgroup~$\Gamma$, but the bound may depend on~$\Gamma$, and the unipotent elementary matrices may not generate quite all of~$\Gamma$:

\begin{thm}[(Carter-Keller-Paige)] \label{CKPThm}
If\/ $\Gamma$ is a finite-index subgroup of\/ $\SL(3,\integer)$, then there is a number~$N$, and a finite-index subgroup\/~$\Gamma'$ of\/~$\Gamma$, such that every element of\/~$\Gamma'$ is a product of\/ $\le N$ unipotent elementary matrices in\/~$\Gamma'$.
\end{thm}

Combining \cref{ElemOrbsBdd,CKPThm} yields the following conclusion:

\begin{cor} \label{GammaHasFPonR}
If\/ $\Gamma$ is a finite-index subgroup of\/ $\SL(3,\integer)$, then every orientation-preserving action of\/~$\Gamma$ on\/~$\real$ has a fixed point.
\end{cor}

\begin{proof}
It suffices to show that the orbit $\Gamma \cdot 0$ is bounded above, for then the supremum of this orbit is a fixed point.

For simplicity, let us ignore the difference between $\Gamma'$ and~$\Gamma$ in \Cref{CKPThm}, so there is a sequence $g_1,g_2,\ldots,g_M$ of unipotent elementary matrices in~$\Gamma$, such that 
	$$ \Gamma = \langle g_M \rangle \, \cdots \, \langle g_2 \rangle \, \langle g_1 \rangle .$$

We will show, for $k = 1,2,\ldots,M$, that
	$$  \text{$\bigl( \langle g_k \rangle \, \cdots \, \langle g_2 \rangle \, \langle g_1 \rangle \bigr) \cdot 0$ is bounded above.} $$
To this end, let $x$ be the supremum of $\bigl( \langle g_{k-1} \rangle \, \cdots \, \langle g_2 \rangle \, \langle g_1 \rangle \bigr) \cdot 0$, and assume, by induction, that $x < \infty$. Then $x \in \real$, so \cref{ElemOrbsBdd} tells us that $\langle g_k \rangle \cdot x$ is bounded above by some $y \in \real$. Then all of $\bigl( \langle g_k \rangle \, \cdots \, \langle g_2 \rangle \, \langle g_1 \rangle \bigr) \cdot 0$ is bounded above by~$y$, as desired.
\end{proof}

\begin{proof}[\normalfont\textbf{Proof of \cref{Witte-AOCThm}}]
Suppose $\Gamma$ has a nontrivial, orientation-preserving action on~$\real$. The set of fixed points is closed, so its complement is a disjoint union of open intervals; let $I$ be one of those intervals.

Note that the interval $I$ is $\Gamma$-invariant (because the endpoints~$I$ are fixed points), so $\Gamma$ acts on~$I$ (by orientation-preserving homeomorphisms).
Then, since the open interval~$I$ is homeomorphic to~$\real$, \cref{GammaHasFPonR} tells us that $\Gamma$ has a fixed point in~$I$.
This contradicts the fact that, by definition, $I$ is contained in the complement of the fixed point set.
\end{proof}

\begin{rem} \label{MayProveNonCpct}
It is hoped that the approach described in this section will (soon?) yield a proof of \cref{LattSLnRNoActConj} in the case where the lattice $\Gamma$ is \emph{not} cocompact.

More precisely, let $\Gamma$ be a lattice in $G = \SL(3,\real)$, such that $G/\Gamma$ is not compact, and suppose we have an orientation-preserving action of~$\Gamma$ on~$\real$. Then:
	\noprelistbreak
	\begin{enumerate}
	\item A matrix~$u$ in $\SL(3,\real)$ is \defit{unipotent} if it is conjugate to an element of $\begin{bmatrix} 1&*&*\\ &1&*\\ &&1 \end{bmatrix}$.
	\item \label{MayProveNonCpct-BddOrb}
\Cref{ElemOrbsBdd} can be generalized to show that if $u$~is any unipotent matrix in~$\Gamma$, then the $u$-orbit of each point in~$\real$ is a bounded set.
	\item \label{MayProveNonCpct-GddGen}
It is well known that some finite-index subgroup~$\Gamma'$ of~$\Gamma$ is generated by unipotent matrices, and it is conjectured that \cref{CKPThm} generalizes: every element of~$\Gamma'$ should be the product of a \emph{bounded} number of unipotent matrices in~$\Gamma'$. 
	\end{enumerate}
If the conjecture in \pref{MayProveNonCpct-GddGen} can be proved, then the above argument shows that $\Gamma$ has no nontrivial, orientation-preserving action on~$\real$.
\end{rem}

\begin{notes}
This combination of bounded orbits and bounded generation was used in \cite{LifschitzMorris-nonarch,LifschitzMorris-BddGen} to prove that some lattices cannot act on~$\real$.

\Cref{SL3ZBddGen} is due to D.\,Carter and G.\,Keller \cite{CarterKeller-BddElemGen}; an elementary proof is given in \cite{CarterKeller-ElemExp}. \Cref{CKPThm} is due to D.\,Carter, G.\,Keller, and E.\,Paige \cite{CarterKellerPaige,Morris-CKP}.

See the remarks leading up to \cref{ReducetoSL3} for more discussion along the lines of \cref{MayProveNonCpct}.

\end{notes}

\section{Complements} \label{ComplementsSect}

\subsection{Actions of lattices in other semisimple Lie groups}

\Cref{LattSLnRNoActConj} refers only to lattices in $\SL(n,\real)$.
We can replace $\SL(n,\real)$ with any other (connected, linear) simple Lie group~$G$ whose real rank is at least~$2$, but some care is needed in stating a precise conjecture for groups that are semisimple, rather than simple. First of all, it should be assumed that the lattice~$\Gamma$ is \defit{irreducible} (i.e., that no finite-index subgroup of~$\Gamma$ is a direct product $\Gamma_1 \times \Gamma_2$ of two infinite subgroups). But additional care is needed if $\SL(2,\real)$ is one of the simple factors of~$G$:

\begin{eg} \label{LattSL2xSL2Acts}
Let $G = \SL(2,\real) \times \SL(2,\real)$, so $G$ is a semisimple Lie group, and $\Rrank G \ge 2$.
Since $\SL(2,\real)$ acts on $\real \cup \{\infty\} \iso \circle$ (by linear-fractional transformations), 
the group~$G$ also acts on~$\circle$, via projection to the first factor. It is then easy to see that any lattice~$\Gamma$ in~$G$ has an action on~$\circle$ (by linear-fractional transformations) in which \emph{every} orbit is infinite. Furthermore, the action is faithful if $\Gamma$ is torsion free and irreducible.
\end{eg}

\begin{conjecture} \label{HighRankFinOrbConj}
Let\/ $\Gamma$ be an irreducible lattice in a connected, semisimple Lie group~$G$ with finite center, such that $\Rrank G \ge 2$.
Then:
	\noprelistbreak
	\begin{enumerate}
	\item $\Gamma$ has no nontrivial, orientation-preserving action on\/~$\real$,
	and
	\item $\Gamma$ is not left orderable.
	\end{enumerate}
Furthermore, if no simple factor of~$G$ is locally isomorphic to $\SL(2,\real)$, then:	\noprelistbreak
\begin{enumerate} \setcounter{enumi}{2}
	\item $\Gamma$ has no faithful action on~$\circle$, 
and
	\item whenever\/~$\Gamma$ acts on~$\circle$, every orbit is finite.
	\end{enumerate}
\end{conjecture}

\begin{rem} \label{AOCConjTrueCasesRem}
The conjecture has been verified in some cases:
	\noprelistbreak
	\begin{enumerate} \renewcommand{\theenumi}{\Alph{enumi}}
	\item \label{AOCConjTrueCasesRem-SL3Sp4}
 \Cref{Witte-AOCThm} verifies the conjecture in the special case where $\Gamma$~is a finite-index subgroup of $\SL(n,\integer)$, with $n \ge 3$. A very similar argument applies when $\Gamma$ is a finite-index subgroup of $\Sp(2n,\integer)$, with $n \ge 2$.
	\item \label{AOCConjTrueCasesRem-Qrank}
From the examples in \pref{AOCConjTrueCasesRem-SL3Sp4}, it follows that the conclusions of the conjecture hold whenever $\Qrank \Gamma \ge 2$.
	\item \label{AOCConjTrueCasesRem-SL2Factor}
L.\,Lifschitz and D.\,W\,Morris verified the conjecture in the special case where 
		\noprelistbreak
		\begin{enumerate}
		\item some simple factor of~$G$ is locally isomorphic to either $\SL(2,\real)$ or $\SL(2,\complex)$,
		and
		\item $G/\Gamma$ is \emph{not} compact.
		\end{enumerate}
	\end{enumerate}
\end{rem}

It seems likely that the method of \cref{BddOrbBddGenSect} will be able to prove the following cases of the conjecture:

\begin{conjecture} \label{SL3RSL3C}
If\/ $\Gamma$ is any non-cocompact lattice in either\/ $\SL(3,\real)$ or\/ $\SL(3,\complex)$, then\/ $\Gamma$ has no nontrivial, orientation-preserving action on\/~$\real$.
\end{conjecture}

If so, then we would have a proof of all of the non-cocompact cases:

\begin{thm}[(Chernousov-Lifschitz-Morris)] \label{ReducetoSL3}
Assume \cref{SL3RSL3C} is true. If $G$ and~$\Gamma$ are as in \cref{HighRankFinOrbConj}, and $G/\Gamma$ is not compact, then the conclusions of \cref{HighRankFinOrbConj} are true.
\end{thm}

Other evidence for \cref{HighRankFinOrbConj} is provided by the Ghys-Burger-Monod Theorem \pref{GhysBurgerMonodFinOrb}, which remains valid in this setting:

\begin{thm}[(Ghys, Burger-Monod)] \label{GhysBurgerMonod-general}
If $G$ and\/~$\Gamma$ are as in \cref{HighRankFinOrbConj}, and no simple factor of~$G$ is locally isomorphic to\/ $\SL(2,\real)$, then:
	\noprelistbreak
	\begin{enumerate}
	\item Every action of\/~$\Gamma$ on~$\circle$ has at least one finite orbit.
	\item $\Gamma$ has no faithful $C^1$ action on the circle~$\circle$.
	\end{enumerate}
\end{thm}

\begin{rem} \label{GhysBurgerMonodGeneralRem} \ 
	\noprelistbreak
	\begin{enumerate}
	\item The lattices that appear in  \cref{HighRankFinOrbConj} are examples of arithmetic groups, and the conjecture could be extended to the class of $S$-arithmetic groups. 
	The Ghys-Burger-Monod Theorem has been generalized to this setting \cite[Cor.~6.11]{WitteZimmer-ActOnCircle}, 
and the appropriate analogue of \cref{HighRankFinOrbConj} has been proved in the special case of $S$-arithmetic groups that are neither arithmetic nor cocompact \cite{LifschitzMorris-nonarch}.
	\item The Burger-Monod proof of \cref{GhysBurgerMonod-general} applies to some cases where $\Gamma$ is a lattice in a product group $G = G_1 \times G_2$ that is not assumed to be a Lie group. A generalization of this result on lattices in products has been proved by U.\,Bader, A.\,Furman, and A.\,Shaker \cite{BaderFurmanShaker}.
	\item \label{GhysBurgerMonodGeneralRem-criterium}
The Burger-Monod injectivity theorem \cf{BurgerMonodInjectSLn} does not immediately imply the conclusion of \cref{GhysBurgerMonod-general} in cases where $H^2(\Gamma;\real) \neq 0$. However, elsewhere in this volume, M.\,Burger \cite{Burger-criterium} uses the injectivity to obtain a general theorem that includes both \cref{BurgerMonodInjectSLn} and the results on lattices in products mentioned in the preceding paragraph. The rough idea is that if we have a $\Gamma$-action on~$\circle$, such that the real Euler class of the action is in the image of $\Hbdd^2(G;\real)$, then a certain quotient of the $\Gamma$-action must extend to a nontrivial, continuous action of~$G$ on the circle.
	\end{enumerate}
\end{rem}

\begin{rem}
Our discussion deals only with actions of lattices on the $1$-dimensional manifolds $\circle$ and~$\real$, but it is conjectured that large lattices also have no faithful actions on manifolds of other small dimensions. (For example, if $n \ge m + 2$, then no lattice in $\SL(n,\real)$ should have a faithful, $C^\infty$ action on any compact $m$-manifold.)
Some discussion of this can be found in 
D.\,Fisher's survey paper \cite{Fisher-Survey} in this volume, or in Robert J.\ Zimmer's CBMS lectures \cite{ZimmerCBMS}.
\end{rem}

\begin{subsecnotes}
A weaker version of \cref{HighRankFinOrbConj} was suggested by D.\,Witte in 1990 (unpublished), but the definitive statement that deals correctly with $\SL(2,\real)$ factors is due to \'E.\,Ghys \cite[p.~200]{GhysCircle}. 

Parts \ref{AOCConjTrueCasesRem-SL3Sp4} and~\ref{AOCConjTrueCasesRem-Qrank} of \cref{AOCConjTrueCasesRem} are due to D.\,Witte \cite{WitteCircle}. See \cite{LifschitzMorris-nonarch} or \cite{LifschitzMorris-BddGen} for Part~\ref{AOCConjTrueCasesRem-SL2Factor}.

\Cref{ReducetoSL3} is implicit in \cite[\S8]{LifschitzMorris-BddGen}. The proof depends crucially on the main result of \cite{ChernousovLifschitzMorris}.

\Cref{GhysBurgerMonod-general} is due to \'E.~Ghys
\cite{GhysCircle} and (in slightly less generality) M.\,Burger and N.\,Monod \cite{BurgerMonod-BddCoho1,BurgerMonod-BddCoho2}. 
Ghys's proof takes a geometric approach that relies on a case-by-case analysis of the possible Lie groups~$G$; a modified (more algebraic) version of the proof that eliminates the case-by-case analysis was found by D.\,Witte and R.\,J.\,Zimmer
\cite{WitteZimmer-ActOnCircle}.
\end{subsecnotes}

\subsection{Some lattices that do act on the circle}

The  following two conjectures suggest that the conclusions of
\cref{HighRankFinOrbConj} fail for lattices in $\SO(1,n)$. Thus,
the assumption that $\Rrank G \ge 2$ cannot be omitted,
although it may be possible to weaken it.

\begin{conjecture}[(W.\,Thurston)] \label{SO1nBettiConj}
 If\/ $\Gamma$ is any lattice in\/ $\SO(1,n)$, then there are 
 \noprelistbreak
 \begin{itemize}
 \item a finite-index subgroup $\Gamma'$ of~$\Gamma$,
 and
 \item a surjective homomorphism $\phi \colon \Gamma' \to \integer$.
 \end{itemize}
 \end{conjecture}

Because $\integer$ obviously has a $C^\infty$~action on the circle with no finite orbits, this implies the following conjecture:

\begin{conjecture} \label{SO1nActConj}
 If\/ $\Gamma$ is any lattice in\/ $\SO(1,n)$, then some finite-index subgroup of\/~$\Gamma$ has a $C^\infty$ action on~$\circle$ that has no finite orbits.
 \end{conjecture}

These conjectures have been proved almost completely, under the additional
assumption that $\Gamma$ is arithmetic.

\begin{thm}[(Li, Millson, Raghunathan, Venkataramana)] \label{SO1nBettiThm}
Suppose\/ $\Gamma$ is a lattice in\/ $\SO(1,n)$. If 
 \noprelistbreak
 \begin{itemize}
 \item $\Gamma$ is arithmetic,
 and
 \item $n \notin \{1,3,7\}$, 
 \end{itemize}
 then the conclusions of \cref{SO1nBettiConj,SO1nActConj} hold.
 \end{thm}

\begin{rem} \label{SO1nBettiRems} \ 
 \noprelistbreak
 \begin{enumerate}
 \item \label{SO1nBettiRems-SO13}
 There exist lattices in $\SO(1,3)$ that act \emph{faithfully}
on~$\circle$. 
(In fact, any torsion-free, cocompact lattice in $\SO(1,3)$ with infinite abelianization is left orderable.)
It would be very interesting to know whether or not there
exist such lattices in $\SO(1,n)$ for $n \ge 4$, or in other groups
of real rank one. 
 \item \label{SO1nBettiRems-hybrid}
 The conclusions of \cref{SO1nBettiConj,SO1nActConj}
hold for the non-arithmetic lattices constructed by M.\,Gromov and I.\,PiatetskiShapiro\cite{GromovPiatetskiShapiro}. Thus, a counterexample to these conjectures
would have to be constructed by some other method.
 \end{enumerate}
 \end{rem}

\begin{rem} \label{SL2Factor->act}
Generalizing \cref{LattSL2xSL2Acts}, it is clear that if $\Gamma$ is a torsion-free, irreducible lattice in~$G$, and $G$ has a simple factor that is locally isomorphic to $\SL(2,\real)$, then $\Gamma$ has a faithful action on~$\circle$ by linear-fractional transformations.
Conversely, under the additional assumption that $\Rrank G \ge 2$, \'E.\,Ghys \cite[Thm.~1.2]{GhysCircle} proved that every action of~$\Gamma$ on~$\circle$ either has a finite orbit or is semi-conjugate to such an action by linear-fractional transformations.
 \end{rem}

\begin{subsecnotes}

\Cref{SO1nBettiConj} is attributed to W.~Thurston
(see \cite[p.~88]{Borel-CohoAsterisque}). 

\Cref{SO1nBettiThm}
combines work of several authors  \cite{Millson,JSLie,LiMillson,RaghunathanVenkataramana}. See \cite{RaghuBettiSurvey} for a survey.

See \cite[\S7.4]{GhysCircleSurvey} for a discussion of some lattices in $\SO(1,3)$ that act on the circle. The general fact stated in \fullcref{SO1nBettiRems}{SO13} is a theorem of S.\,Boyer, D.\,Rolfsen, and B.\,Wiest \cite{BoyerRolfsenWiest}. 

Remark~\fullref{SO1nBettiRems}{hybrid}
is due to A.~Lubotzky \cite{Lubotzky-hybridBetti}.
\end{subsecnotes}

\subsection{Actions on trees}

The real line and the circle are the only connected $1$-manifolds, but there are many other $1$-dimensional simplicial complexes.
For short, a contractible, $1$-dimensional simplicial  complex is called a \emph{tree}, and focusing our attention on the
groups that act on trees leads to an
interesting theory (the \defit{Bass-Serre theory of group actions on
trees}).

\begin{rem}
Suppose, as usual, that
$\Gamma$ is a lattice in a (connected, linear) semisimple Lie group~$G$.
Then it is easy to construct a faithful action of~$\Gamma$ on some tree, even if we assume that the tree is locally finite. To do this:
	\noprelistbreak
	\begin{enumerate}
	\item Let $\Gamma = N_0 \supset N_1 \supset \cdots$ be a chain of
finite-index, normal subgroups of~$\Gamma$, such that $\bigcap_k N_k =
\{e\}$.
	\item Let the $0$-skeleton $T_0$ be the disjoint union of all $\Gamma/N_k$.
	\item Let the $1$-skeleton $T_1$ have a $1$-simplex (or ``edge'') joining $\gamma N_k$ and~$\gamma N_{k+1}$, for every $\gamma \in \Gamma$ and $k = 0,1,\ldots$.
	\end{enumerate}
Then $\Gamma$ has a natural action on~$T_0$ by left translations, and this extends to an action on~$T$ (by isometries).
\end{rem}

The following theorem states that, under mild hypotheses, every action of~$\Gamma$ on a tree has a finite orbit. This can be thought of as an analogue of \cref{GhysBurgerMonodFinOrb} for actions on trees.

\begin{thm} \label{GammaHasFA}
 If 
 \noprelistbreak
 \begin{itemize}
 \item $\Gamma$ is as in \cref{HighRankFinOrbConj}, or\/ $\Gamma$ has Kazhdan's property~$(T)$,
 \item $T$ is a tree that is not homeomorphic to~$\real$,
 and
 \item $\Gamma$ acts on~$T$ by homeomorphisms,
 \end{itemize}
 then $\Gamma$ has at least one finite orbit on~$T$.
 \end{thm}
 
 \begin{rem} \ 
 \noprelistbreak
 \begin{enumerate}
\item More precisely, the finite orbit in the conclusion of \cref{GammaHasFA} can be taken to consist of either a single vertex or two vertices of the tree.
\item In the Bass-Serre theory, it is usually assumed that the action is by isometries. In this case:
	\noprelistbreak
	\begin{enumerate}
	\item  there is no need to assume that $T$ is not homeomorphic to~$\real$, 
	and 
	\item the finite orbit can be taken to be a fixed point (and this fixed point is either a vertex or the midpoint of some edge).
	\end{enumerate}
\end{enumerate}
 \end{rem}

\begin{rem} \label{FA<>}
 A fundamental conclusion of the Bass-Serre theory is that there is a finite orbit in every action of a countable group~$\Lambda$ on every tree (except possibly~$\real$) if and only if
 \noprelistbreak
 \begin{enumerate}
 \item $\Lambda$ is finitely generated,
 \item $\Lambda/[\Lambda,\Lambda]$ is finite,
 and
 \item $\Lambda$ cannot be written in any nontrivial way as a free
product with amalgamation $A \mathbin{*_C} B$.
 \end{enumerate}
 In the situation of \cref{GammaHasFA}, it is well known that $\Gamma$~is finitely generated, and that $\Gamma/[\Gamma,\Gamma]$ is finite.
Thus, in algebraic terms, \cref{GammaHasFA} is the assertion that
$\Gamma$~is not a free product with amalgamation.
 \end{rem}

 \begin{subsecnotes}
J.--P.\,Serre's elegant book \cite{Serre-Trees} is the standard
introduction to the Bass-Serre theory of actions on trees. 

In the special case where $\Gamma = \SL(3,\integer)$, 
\cref{GammaHasFA} is due to J.--P.~Serre \cite{Serre-SL3ZhasFA}. (See \cite[Thm.~16, p.~67]{Serre-Trees} for an exposition.) The generalization to other lattices of higher rank is due to G.~A.~Margulis \cite[Thm.~2]{MargulisFA}. The case of groups with Kazhdan's property $(T)$ is due to R.~Alperin \cite{Alperin-T->FA}
and Y.~Watatani \cite{Watatani-T->FA}, independently. (Proofs can also
be found in \cite[\S2.3]{BekkaEtAlT} and
\cite[Thm.~3.3.9 and \P3.3.10]{MargulisBook}.)

See \cite[Thm.~15, p.~58]{Serre-Trees} for a proof of \cref{FA<>}.
\end{subsecnotes}

{\scshape Department of Mathematics and Computer Science,
University of Lethbridge, Lethbridge, Alberta, T1K~3M4, Canada}

\end{document}